\documentclass[12pt]{amsart}
\usepackage{amsfonts,amssymb,amscd,amstext,mathrsfs}
\usepackage[utf8]{inputenc}
\usepackage{hyperref}
\usepackage{verbatim}
\usepackage{color} 
\usepackage{times}
\usepackage{enumerate}
\usepackage[up,bf]{caption}
\usepackage{graphicx}
\usepackage[usenames,dvipsnames,svgnames,table]{xcolor}
\usepackage{tikz-cd}
\usepackage{mathrsfs}

\usepackage[top=1.8in, bottom=1.6in, left=1.6in, right=1.6in]{geometry}

\input xy
\xyoption{all}

\newtheorem{theorem}{Theorem}[section]
\newtheorem{proposition}[theorem]{Proposition}

\newtheorem{corollary}[theorem]{Corollary}

\theoremstyle{definition}
\newtheorem{definition}[theorem]{Definition}
\newtheorem{remark}[theorem]{Remark}

\newtheorem{problem}[theorem]{Problem}
\newtheorem{example}[theorem]{Example}

\numberwithin{equation}{section}
\numberwithin{figure}{section}

\newcommand\Ocal{\mathcal{O}}

\newcommand\Scal{\mathcal{S}}

\newcommand{\Hom}{\text{Hom}}
\newcommand\Cscr{\mathscr{C}}

\newcommand\Oscr{\mathscr{O}}

\newcommand\A{\mathbb{A}}
\newcommand\B{\mathbb{B}}
\newcommand\C{\mathbb{C}}
\newcommand\D{\overline{\mathbb D}}
\newcommand\CP{\mathbb{CP}}

\renewcommand\D{\mathbb D}

\newcommand\N{\mathbb{N}}

\newcommand\R{\mathbb{R}}

\newcommand\Z{\mathbb{Z}}
\newcommand{\VF}{ \operatorname{{\rm VF}}}
\newcommand{ \LieO}{ \operatorname{{\rm Lie}_{hol}^\omega}}
\newcommand\ggot{\mathfrak{g}}

\newcommand\igot{\mathfrak{i}}

\renewcommand\igot{\mathfrak{i}}

\newcommand\E{\mathrm{e}}

\renewcommand\imath{\igot}

\newcommand\hra{\hookrightarrow}

\newcommand\wt{\widetilde}
\newcommand\wh{\widehat}
\newcommand\di{\partial}

\renewcommand\div{\mathrm{div}}

\newcommand\Aut{\mathrm{Aut}}
\newcommand\Autalg{\mathrm{Aut}_{\mathrm{alg}}}
\newcommand\Aff{\mathrm{Aff}}
\newcommand\GL{\mathrm{GL}}
\newcommand\SL{\mathrm{SL}}

\newcommand\Id{\mathrm{Id}}

\newcommand\End{\mathrm{End}}

\newcommand\VFH{\mathrm{VF_{hol}}}
\newcommand{\VFalg}[0]{\ensuremath{\operatorname{VF}_{\mathrm{alg}}}}
\newcommand{\Liealg}[0]{\ensuremath{\operatorname{Lie}_{\mathrm{alg}}}}
\newcommand{\Xsing}[0]{\ensuremath{X_{\mathrm{sing}}}}
\newcommand{\TT}[0]{\ensuremath{\mathrm{T}}}
\numberwithin{equation}{section}
\newcommand{\mat}[3]{{\mathrm{Mat}\left( #1 \times #2 ; \, #3 \right)}} %space of matrices

\begin{document}
\title{The first thirty years of Anders\'en--Lempert theory}
\author{Franc Forstneri{\v c} and Frank Kutzschebauch}

\address{Franc Forstneri\v c: Faculty of Mathematics and Physics, University of Ljubljana, and Institute of Mathematics, Physics, and Mechanics, Jadranska 19, 1000 Ljubljana, Slovenia}
\email{franc.forstneric@fmf.uni-lj.si}

\address{Frank Kutzschebauch: Institute of Mathematics, University of Bern, Sidlerstrasse 5, CH-3012 Bern, Switzerland}
\email{Frank.Kutzschebauch@math.unibe.ch}  

\thanks{Forstneri\v c is supported by research program P1-0291 and grant J1-3005 from ARRS, 
Republic of Slovenia.} 
\thanks{Kutzschebauch is supported by Schweizerische Nationalfonds Grant Nr.\ 200021-178730.}


\date{\today}
%%%%%


\begin{abstract}
In this paper we expose the impact of the fundamental discovery, made by 
Erik Anders\'en and L\'aszl\'o Lempert in 1992, that the group generated by shears is dense in the group 
of holomorphic automorphisms of a complex Euclidean space of dimensions $n>1$. 
In three decades since its publication, their groundbreaking work led to the discovery of several 
new phenomena and to major new results in complex analysis and geometry involving Stein manifolds 
and affine algebraic manifolds with many automorphisms. The aim of this survey is to present the 
focal points of these developments, with a view towards the future.
\end{abstract}

\maketitle

\centerline{\em Dedicated to L\'aszl\'o Lempert in honour of his 70th birthday}

\tableofcontents

%
%
%    INTRODUCTION
%
%
\section{Introduction}\label{sec:intro} 
Complex Euclidean spaces $\C^n$ $(n\in \N=\{1,2,3,\ldots\})$ are the most basic and important 
objects in analytic and algebraic geometry. It is natural to try understanding holomorphic 
automorphisms (symmetries) of $\C^n$ and their role in applications. If $n=1$, 
these are precisely the affine linear maps $z\mapsto az+b$ with $a\in \C^*=\C\setminus\{0\}$
and $b\in \C$. However, for $n>1$ the group $\Aut(\C^n)$ of holomorphic automorphisms
of $\C^n$ is huge. The affine group $\Aff(\C^n)$, generated by
the general linear group $\GL_n(\C)$ together with translations, acts transitively on $\C^n$.
Writing complex coordinates on $\C^n$ as $z=(z',z_n)$ with $z'=(z_1,\ldots,z_{n-1})\in \C^{n-1}$, 
we have automorphisms of the form
\begin{equation}\label{eq:shear}
	\Phi(z)=\Phi(z',z_n)=\left(z', \E^{f(z')} z_n+g(z')\right), \quad\ z\in\C^n, 
\end{equation}
where $f$ and $g$ are entire functions of $n-1$ variables. Such maps and their $\GL_n(\C)$-conjugates
are called {\em shears}. Maps \eqref{eq:shear} with $f=0$ and their $\SL_n(\C)$-conjugates
are {\em additive shears}; they have complex Jacobian determinant identically equal to $1$. 
Maps \eqref{eq:shear} with $g=0$ and their $\GL_n(\C)$-conjugates are {\em multiplicative shears}.
Shears generate the {\em shear group} $\Scal(n)$, an infinite dimensional
subgroup of $\Aut(\C^n)$ when $n>1$. Likewise, additive shears generate a 
subgroup $\Scal_1(n)$ of the group $\Aut_1(\C^n)$ of holomorphic automorphisms with Jacobian one.
By composing shears, one obtains many interesting automorphisms.
For example, composing a shear in two variables and the switch map 
$
	\delta(z_1,z_2)=(z_2,z_1)
$
yields {\em H\'enon maps} 
$
	H(z_1,z_2) = \left(\E^{f(z_1)}z_2+g(z_1),z_1\right)\! 
$
whose dynamical properties have been studied intensively.

In 1990, Erik Anders\'en \cite{Andersen1990} made a fundamental discovery that 
for every $n>1$ the group $\Scal_1(n)$ generated by additive shears is dense in 
$\Aut_1(\C^n)$ but not equal to it. This was extended by Anders\'en
and L\'aszl\'o Lempert \cite{AndersenLempert1992} in 1992 to the pair of groups
$\Scal(n)\subset \Aut(\C^n)$. A year later, Forstneri\v c and Rosay \cite{ForstnericRosay1993}
recast the proof of their approximation theorem \cite[Theorem 1.3]{AndersenLempert1992} 
in terms of complete holomorphic vector fields, showing that
the main technical lemma from \cite{Andersen1990,AndersenLempert1992} 
amounts to the following statement: 

\smallskip \label{AL-observation}
(*) {\em Every polynomial holomorphic vector field $V$ on $\C^n$ for $n>1$ is a finite sum of complete 
polynomial vector fields $V_1,\ldots, V_N$ whose flows consist of shears.
If $V$ has divergence zero then each $V_j$ can be chosen to have divergence zero, 
and we only get additive shears.}
\smallskip

The flow of a complete holomorphic vector field is a 
one-parameter group of holomorphic automorphisms; conversely, the infinitesimal generator 
of a one-para\-meter group of automorphisms is a complete holomorphic vector field. 
Since the flow of the sum $V=V_1+V_2+\cdots+V_N$ can be approximated
by compositions of flows of vector fields $V_1,\ldots, V_N$, it follows that the
flow of any holomorphic vector field on $\C^n$ can be approximated
%, uniformly on compacts on its domain of existence, 
by compositions of shears. This implies that the shear group
$\Scal(n)$ is dense in the automorphism group $\Aut(\C^n)$. Forstneri\v c and Rosay gave the following 
more general and useful version of this approximation 
result; see \cite[Theorem 1.1 and Erratum]{ForstnericRosay1993}, 
as well as \cite[Theorem 1.1]{Forstneric1994MA}.

%
%   AL THEOREM
%
\begin{theorem}\label{th:FR} 
Assume that $\Omega$ is a domain in $\C^n$, $n>1$, and $\Phi_t:\Omega\to \C^n$ $(t\in [0,1])$ 
is a continuous isotopy of injective holomorphic maps such that 
$\Phi_0$ is the identity map on $\Omega$ and the domain
$\Omega_t=\Phi_t(\Omega)$ is Runge in $\C^n$ for every $t\in[0,1]$. Then,
$\Phi_1$ can be approximated uniformly on compacts in $\Omega$ by 
elements of the shear group $\Scal(n)$. If in addition the domain $\Omega$ is pseudoconvex, 
$H^{n-1}(\Omega,\C)=0$, and $\Phi_t$ has Jacobian one for every $t\in [0,1]$,
then $\Phi_1$ can be approximated by elements of the shear group $\Scal_1(n)$. 
\end{theorem}

Recall that a domain $\Omega$ in $\C^n$ is called Runge if polynomials are dense in the 
space $\Oscr(\Omega)$ of holomorphic functions. 
The domain $\Omega$ in the first part of Theorem \ref{th:FR} need not be pseudoconvex; 
however, the Runge property is essential. It is easily seen that if $\Phi:\Omega\to \Omega'$ 
is a biholomorphic map between domains in $\C^n$ which is a limit of automorphisms
of $\C^n$, then $\Omega$ is Runge if and only if $\Omega'$ is Runge.
The condition $H^{n-1}(\Omega,\C)=0$ and pseudoconvexity of $\Omega$
guarantee that holomorphic vector fields with vanishing divergence on $\Omega$ 
can be approximated by holomorphic vector fields of the same type defined on $\C^n$. 
(The point is that, on a Stein manifold, the de Rham cohomology can be computed 
by means of holomorphic differential forms.)
A complete proof of Theorem \ref{th:FR} is also available in \cite[Theorem 4.9.2]{Forstneric2017E}. 

The next important step was made by Dror Varolin in his dissertation 
(University of Wisconsin-Madison, 1997) and in the papers 
\cite{Varolin1999,Varolin2000,Varolin2001}. His vantage point is the observation
that the flow of a holomorphic vector field on a complex manifold $X$,
which is a Lie combination % (using sums and Lie brackets)
of complete holomorphic vector fields, is a limit of holomorphic automorphisms 
of $X$. Mimicking (*), Varolin introduced the following notion.

%
%   DENSITY PROPERTY
%
\begin{definition}\label{def:density}
A complex manifold $X$ has the {\em density property} if every holomorphic 
vector field on $X$ can be approximated uniformly on compacts by Lie combinations
of complete holomorphic vector fields on $X$.
\end{definition}

Varolin \cite{Varolin2001} introduced the notion of density property for any Lie algebra of 
holomorphic vector fields (an {\em algebraic structure} on $X$), 
asking that it be densely generated by the complete vector fields
which it contains; see Definition \ref{def:densitygeneral}. An important example is the Lie algebra of
holomorphic vector fields having zero divergence with respect to a holomorphic volume form 
$\omega$ on $X$, that is, a nowhere vanishing holomorphic section of the 
canonical bundle $K_X=\wedge^n T^*X$ with $n=\dim X$ (the top exterior power of the cotangent 
bundle of $X$). The density property for this Lie algebra is called the {\em volume density property} of 
$(X,\omega)$. The standard volume form on $\C^n$ with coordinates $z=(z_1,\ldots,z_n)$ is  
$dz_1\wedge dz_2\wedge\cdots \wedge dz_n$.  % See Section \ref{sec:DP} for more details.

The density property of $X$ is equivalent to the ostensibly weaker condition that 
every holomorphic vector field on $X$ can be approximated uniformly 
on compacts by sums of complete holomorphic vector fields on $X$; however,
commutators come handy in calculations. 
Thus, the Anders\'en--Lempert observation (*) says that $\C^n$ for $n>1$ enjoys the density property.
By mimicking the proof in the case $X=\C^n$, one easily obtains the following analogue
of Theorem \ref{th:FR}; see \cite[Theorem 4.10.5]{Forstneric2017E}.

%
%
%  FORSTNERIC AND ROSAY FOR STEIN MANIFOLDS WITH DENSITY PROPERTY
%
%
\begin{theorem}\label{th:FR2}
Let $X$ be a Stein manifold with the density property. If 
$\Phi_t : \Omega \stackrel{\cong}{\to} \Omega_t=\Phi_t(\Omega)\subset X$ $(t\in [0,1])$ 
is a continuous isotopy of biholomorphic maps between Stein  
Runge domains in $X$ with $\Phi_0=\Id_{\Omega}$,  then $\Phi_1$ can be approximated 
uniformly on compacts in $\Omega$ by holomorphic automorphisms of $X$.
\end{theorem}

The assumption that the domains $\Omega_t$ are Runge and Stein  
is used to approximate holomorphic vector fields on $\Omega_t$ by holomorphic vector fields on $X$. 
If $X=\C^n$ or, more generally, if $X$ is holomorphically parallelizable, this reduces to 
approximation of functions, and hence it suffices to assume that $\Omega_t$ 
is Runge in $X$ for each $t\in [0,1]$.

%
%  REMARK:  COMPLETNESS IN REAL TIME
%
\begin{remark}
In this paper, a holomorphic vector field $V$ on a complex manifold $X$ is said to be 
complete it is complete in complex time, i.e., its flow is a complex $1$-parameter group
of holomorphic automorphisms of $X$. One could also consider the 
density property for holomorphic vector fields which are complete
in real time. However, it turns out that every Stein manifold $X$ with this {\em real density
property} also has the density property as defined above. Indeed,
working with $\R$-complete holomorphic vector fields, the proof of Theorem \ref{th:FR2} 
holds without any changes. Therefore, we can find a Fatou--Bieberbach domain in $X$
around each point $x\in X$, arising by contracting a small coordinate ball around $x$ and 
approximating this map by an automorphism of $X$ with an attracting fixed point at $x$. The basin of 
attraction is biholomorphic to $\C^n$ with $n=\dim X$ (see Theorem \ref{th:basin}). It follows that 
the manifold $X$ has the Liouville property, i.e., every negative plurisubharmonic function on $X$ is 
constant. On a Stein manifold with the Liouville property, every  holomorphic vector field which is 
complete in real time is also complete in complex time, see \cite[Corollary 2.2]{Forstneric1996MZ}. 
\end{remark}

Here is another useful version of Theorems \ref{th:FR} and \ref{th:FR2} for isotopies of compact 
holomorphically convex sets; see \cite[Theorem 2.1]{ForstnericRosay1993} or 
\cite[Theorem 4.12.1]{Forstneric2017E} for the case $X=\C^n$, $n>1$. The proof in
the general case is the same. 

%
%   ISOTOPIES OF COMPACT HC SETS
%
\begin{theorem}\label{th:ALPC}
Let $X$ be a Stein manifold with the density property. Assume that 
$K$ is a compact holomorphically convex (i.e., $\Oscr(X)$-convex) set in $X$,
$\Omega\subset X$ is an open set containing $K$, and $\Phi_t:\Omega\to X$
$(t\in[0,1])$ is a continuous isotopy of injective holomorphic maps such that $\Phi_0=\Id_\Omega$
and the set $K_t=\Phi_t(K)$ is holomorphically convex in $X$ for every $t\in [0,1]$.
Then we can approximate $\Phi_1$ uniformly on $K$ by holomorphic automorphisms of $X$. 
\end{theorem}
 
There is an easy reduction to Theorem \ref{th:FR2} since a compact holomorphically  
convex set admits a basis of Stein Runge neighbourhoods in $X$.
In both version of this result, one can approximate the entire isotopy $\{\Phi_t\}_{t\in[0,1]}$
by isotopies of automorphisms of $X$. 

Several further additions to Theorems \ref{th:FR} and \ref{th:FR2} are possible. 
The following parametric version of Theorem \ref{th:FR} was proved by Kutzschebauch
\cite[Theorem 2.3]{Kutzschebauch2005}.  

%
%   PARAMETRIC VERSION
%
\begin{theorem}\label{th:parameters}
Let $\Omega$ be an open set in $\C^n=\C^k\times \C^m$ with $m>1$.
For every $t \in [0, 1]$ let $\Phi_t:\Omega \to \C^n$ be a continuous family of 
injective holomorphic map of the form
\begin{equation}
\label{eq:aut-par}
        \Phi_t (z,w) = \bigl(z,\varphi_t (z,w)\bigr), \quad z\in \C^k,\ w\in \C^m
\end{equation}
%and of class ${\mathcal C}^1$ in $(t, z, w) \in [0, 1]\times \Omega$, 
with $\Phi_0=\Id_\Omega$. If the domain $\Phi_t (\Omega)$ is Runge in $\C^n$ for every $t\in [0,1]$,
then $\Phi_1$ can be approximated uniformly on compacts in $\Omega$ by holomorphic
automorphisms of the form \eqref{eq:aut-par}.
\end{theorem}

It is elementary to include interpolation at finitely many points in  
these approximation theorems. There are considerably more general interpolation results for automorphisms
on algebraic subvarieties of $\C^n$, which we describe in Sections \ref{sec:DP} and \ref{sec:jets}.
In particular, any jet of a local biholomorphism at a point is the jet of an automorphism;
see Theorem \ref{th:BF2000} and its parametric version, Theorem \ref{th:RPU}.
These results are useful the construction of holomorphic automorphisms
with interesting dynamical properties.

Theorem \ref{th:FR} was extended by Kutzschebauch and Wold 
\cite{KutzschebauchWold2018} (2018) to include Carleman approximation on sets of the form 
$K\cup \R^s\subset \C^n$, where $\R^s$ is an affine totally real subspace of $\C^n$ for $s<n$. 
Under certain technical assumptions on the isotopy of such sets in $\C^n$ (they must be 
polynomially convex, totally real on $\R^s\setminus K$, and nearly fixed near infinity on $\R^s$), 
it is possible to approximate the isotopy in the Carleman sense by isotopies of 
holomorphic automorphisms of $\C^n$. 

Theorems \ref{th:FR} and \ref{th:FR2} and their generalisations mentioned in the sequel 
form the backbone of what is now called {\em the Anders\'en--Lempert theory}.

%
%   MANIFOLDS WITH DENSITY PROPERTY; ELLIPTIC COMPLEX GEOMETRY
%
After the initial work of Varolin, the subject of Stein manifolds with the density 
(or volume density) property was developed by Kaliman and Kutzschebauch with collaborators; 
see Section \ref{sec:DP} for an overview. They also introduced and studied the algebraic 
(volume) density property on affine complex manifolds, the algebraic analogues of Stein manifolds. 
A closely related notion is {\em holomorphically flexibility} in the sense of Arzhantsev et al.\ 
\cite{ArzhantsevFlennerKalimanKutzschebauchZaidenberg2013DMJ}, asking that complete
holomorphic vector fields generate the tangent space of the manifold at every point.
This field became known as {\em elliptic complex geometry}.
Most complex Lie groups and homogeneous manifolds have the density property.
In the two decades since their introduction,
Stein manifolds with the density property came to play a major role in complex geometry. 
Every such manifold is an Oka manifold (see Section \ref{sec:Oka}); however,
the former class of manifolds allows many additional and more precise results concerning 
holomorphic maps from Stein manifolds. On Stein manifolds with the density property 
one can find automorphisms with given jets at finitely many points, and also at 
some discrete set of points; see Section \ref{sec:jets}.

Sections \ref{sec:FB}--\ref{sec:3D} contain a survey of applications of Anders\'en--Lempert theory. 

%
%	FB DOMAINS
%
A fascinating phenomenon in complex analysis
is the existence for any $n>1$ of proper subdomains $\Omega\subsetneq\C^n$ 
which are biholomorphic to $\C^n$; see Section \ref{sec:FB}. 
Such a domain is called a {\em  Fatou--Bieberbach domain}.
An injective holomorphic map $\Phi:\C^n\to \C^n$ whose image $\Omega=\Phi(\C^n)$
is a proper subdomain of $\C^n$ (and its inverse $\Phi^{-1}:\Omega\to\C^n$)  
is called a {\em Fatou--Bieberbach map}. 
The first explicit example was given by Bieberbach \cite{Bieberbach1933} in 1933,
following earlier examples by Fatou \cite{Fatou1922} (1922) of non-degenerate 
(but non-injective) entire maps $\C^2\to\C^2$ with non-dense images.
(See also Bochner and Martin \cite[Sect. III.1]{BochnerMartin1948}.)
These examples arose as limits of holomorphic automorphisms of $\C^n$, and hence
they are Runge in $\C^n$.
%Indeed, if a sequence $\Phi_k\in\Aut(\C^n)$ converges uniformly on compacts in $\C^n$ and the limit map $\Phi=\lim_{k\to\infty}\Phi_k:\C^n\to\C^n$ has maximal rank $n$ at some point, then $\Phi$ is an injective holomorphic map which however need not be surjective when $n>1$. If a Fatou--Bieberbach map $\Phi:\C^n \stackrel{\cong}{\to} \Omega\subsetneq\C^n$ is a limit of automorphisms of $\C^n$ then $\Omega$ is Runge in $\C^n$.
In particular, if $p\in\C^n$ is an attracting fixed point of a holomorphic automorphism
$F\in\Aut(\C^n)$ and we denote by $F^k$ the $k$-th iterate of $F$, the basin of attraction
$
	\{z\in \C^n: \lim_{k\to\infty} F^{k}(z) =p\}
$
is either $\C^n$ or a Runge Fatou--Bieberbach domain (see Theorem \ref{th:basin}).
Many further examples of Fatou--Bieberbach maps with interesting properties arose as limits 
of sequences of automorphisms which are not iterates of a fixed automorphism.
Theorems \ref{th:FR} and \ref{th:FR2} gave rise to several new constructions of 
Fatou--Bieberbach domains. An interesting application of  the 
Anders\'en--Lempert theory was the discovery by Wold \cite{Wold2008} (2008)  
of Fatou--Bieberbach domains in $\C^n$ for any $n>1$ which fail to be Runge, and
hence the corresponding Fatou--Bieberbach maps are not limits of automorphisms of $\C^n$. 

In Section \ref{sec:twisted} we describe constructions of highly twisted proper holomorphic 
embeddings $F:\C^k\hra\C^n$ for $1\le k<n$ such that the complements $\C^n\setminus F(\C^k)$
of their images are $(n-k)$-hyperbolic in the sense of Eisenman 
(see Theorem \ref{th:hyperboliccomplement}). In particular, 
$\C^k$ embeds as a closed complex hypersurface $\Sigma\subset \C^{k+1}$ 
such that $\C^{k+1}\setminus \Sigma$ is Kobayashi hyperbolic.
This complements the well-known fact that most affine algebraic hypersurfaces 
in $\C^n$ of sufficiently big degree are hyperbolic and have hyperbolic complements. 
This led to the discovery by Derksen and Kutzschebauch 
\cite{DerksenKutzschebauch1998} of nonlinearizable periodic holomorphic automorphisms 
of $\C^n$ for any $n\ge 4$; see Theorem \ref{th:nonlinearizable}.

%
%  EMBEDDING OPEN RIEMANN SURFACES IN C^2
%
Another classical subject where the Anders\'en--Lempert theory generated major progress is
the {\em Forster--Bell--Narasimhan Conjecture} \cite{Forster1970,BellNarasimhan1990},
asking whether every open Riemann surface embeds properly as a smooth complex curve in the
Euclidean plane $\C^2$. The most advanced known results on this subject
use constructions of Fatou--Bieberbach domains with the aid of Theorem \ref{th:FR}.
We discuss  this topic in Section \ref{sec:RS}; a more comprehensive presentation 
can be found in \cite[Sections 9.10--9.11]{Forstneric2017E}. 

%
%   LONG C^n's
%
Wold's construction in \cite{Wold2010} of non-Runge Fatou--Bieberbach domains
led to the construction of uncountably many pairwise non-biholomorphic 
and non-Stein long $\C^n$'s, i.e., complex manifolds which are increasing 
unions of biholomorphic copies of $\C^n$; see Wold \cite{Wold2010} and
Boc Thaler and Forstneri\v c \cite{BocThalerForstneric2016}. 
Many interesting question concerning these {\em long $\C^n$'s} remain open,
and we refer to Section \ref{sec:long} for a discussion of this topic.  

In Section \ref{sec:Oka} we discuss the impact of the Anders\'en--Lempert theory on 
Oka theory. Every Stein manifold $X$ with the density or volume density property is an Oka manifold
(Theorem \ref{th:DPOka}), 
which means that holomorphic maps $S\to X$ from any Stein manifold $S$ satisfy
all forms of the h-principle (see \cite[Theorem 5.4.4]{Forstneric2017E}). 
Furthermore, if $K\subset X$ is a compact $\Oscr(X)$-convex subset then $X\setminus K$ is an 
Oka manifold (Theorem \ref{th:KOka}). The same holds for complements of certain closed 
unbounded $\Oscr(X)$-convex subsets. However, on a Stein manifold $X$ 
the density property is a much stronger condition than the Oka property, implying 
finer results on the existence of holomorphic maps $S\to X$ which do not hold for every 
Oka manifold $X$. For example, a Stein manifold $X$ with the density property
contains every Stein manifold of dimension $k$ with $2k+1\le \dim X$ as a properly embedded 
complex submanifold (immersed if $2k=\dim X$), and one can even prescribe the homotopy class
of the embedding. Furthermore, Stein manifolds with the density property 
contain big Stein Runge domains which are total spaces of normal bundles of certain embedded
complex submanifolds. For example, $\C^*\times \C$ embeds as a Runge domain in $\C^2$,
although not in any obvious way. 

In Section \ref{sec:complete} we survey recent results concerning the problem
of Paul Yang from 1977, asking whether there are bounded (metrically) complete complex 
submanifolds of $\C^n$. (This is holomorphic analogue of the Calabi--Yau problem
concerning minimal surfaces in $\R^n$ for $n\ge3$; 
see \cite[Chapter 7]{AlarconForstnericLopez2021} for the latter.)  
It has been discovered fairly recently that the ball of $\C^n$ and, more generally, 
any pseudoconvex Runge domain in $\C^n$ can be foliated by complete complex submanifolds 
of any codimension and with partial control of the topology of the leaves.
The methods of Anders\'en--Lempert theory play a major role in these constructions.

In Section \ref{sec:3D} we discuss an application of the Anders\'en--Lempert method in the
smooth world. A long standing problem in 3-dimensional topology asked whether the fundamental 
group of any homology 3-sphere different from the 3-sphere $S^3$ admits an irreducible 
representation into $\mathrm{SL}_2 (\C)$, i.e., a 2-dimensional irreducible representation. 
%(Recall that a homology 3-sphere is a compact 3-manifold $X$ whose homology groups are those of $S^3$. Its fundamental group is nontrivial unless $X$ is $S^3$.)
The affirmative answer given by Rafael Zentner \cite{Zentner2018} in 2018 is a case where shears 
play a role in the real setting.  

We conclude by discussing the {\em recognition problem in complex analysis} in 
Section \ref{sec:recognition}. The question is how to decide whether a given Stein manifold, 
which is contractible and simply connected at infinity, is a complex Euclidean space. 
To be such, it must have many holomorphic automorphisms. Hence, it is natural to ask
whether every Stein manifold with the density property which is diffeomorphic to $\R^{2n}$ 
is also biholomorphic to $\C^n$ (see Problem \ref{prob:TothVarolin}).
This question, asked by T\'oth and Varolin \cite{TothVarolin2000} in 2000, 
remains unsolved. We also discuss other related problems such as the cancellation problem.

%
%
% SECTION: STEIN MANIFOLDS WITH THE DENSITY PROPERTY
%
%
\section{Stein manifolds with density properties}\label{sec:DP}

We denote by $\VFH (X)$ the Lie algebra of global holomorphic vector fields on a 
complex manifold $X$. In \cite{Varolin2001} Varolin introduced the following notion.

\begin{definition} \label{def:densitygeneral}
A Lie subalgebra  $\ggot$ of the Lie algebra $\VFH (X)$ 
has the density property if the Lie algebra  generated by complete fields 
in $\ggot$ is dense in $\ggot$ in the compact-open topology.
\end{definition}

%
%	SUBSECTION: DENSITY PROPERTY
%
\subsection{Density property}\label{ss:DP}
The case that $\ggot = \VFH (X)$ 
corresponds to $X$ having the density property (see Definition \ref{def:density}). 
The importance of this property for Stein manifolds lies in Theorem \ref{th:FR2} from 
the introduction. On the other hand, for compact complex manifolds 
this property is of no interest since every holomorphic vector field is complete, and
furthermore there need not exist any nonzero holomorphic vector fields.

Establishing the density property can be rather tricky. 
After the initial work of Anders\'en and Lempert \cite{AndersenLempert1992} 
which established this property for Euclidean spaces $\C^n$ of any dimension $n\ge 2$, 
Varolin \cite{Varolin2000,Varolin2001} gave a list of further examples, 
and in his joint work with T\'oth \cite{TothVarolin2000} established the density property 
for semi-simple complex Lie groups and their quotients by reductive subgroups.
Later, Kaliman and Kutzschebauch \cite{KalimanKutzschebauch2008IM}
found a strong criterion based on the notion of compatible pairs of vector fields;
see Definitions \ref{def:compatible-alg} and \ref{def:compatible}.  
Their criterion was initially developed for affine algebraic manifolds.
For such manifolds, the algebraic density property was already introduced by Varolin
as follows.

\begin{definition}\label{def:algebraicdensity}
An affine algebraic manifold $X$ has the algebraic density property if 
the Lie algebra of algebraic vector fields on $X$ is generated by complete algebraic vector fields.
\end{definition}

Varolin remarked that this condition implies the holomorphic density property.
Indeed, the Oka--Weil theorem says that polynomial functions 
are dense in the space of holomorphic functions; 
by an application of Theorems A and B the same holds for sections of any coherent 
algebraic sheaf, in particular, for holomorphic vector fields on an affine algebraic manifold.

%
%  COMPATIBLE PAIRS
%
\begin{definition} \label{def:compatible-alg}
Complete algebraic vector fields $\nu$ and $\mu$ on an affine algebraic manifold $X$ 
form a {\em compatible pair} if the following two conditions hold:
\begin{enumerate}[(1)] 
\item the linear span of the product of the kernels $\ker \nu \,\cdotp \ker \mu$ contains a 
nontrivial ideal $I \subset \C[X]$, and
\item there is a function $h \in \C[X]$ with  $h\in \ker \mu$ and $\nu (h) \in \ker \nu\setminus \{0\}$.  
\end{enumerate}
 If only condition (1) is satisfied, we call $(\nu,\mu)$ a semi-compatible pair.
The biggest ideal $I$ with this property is called the {\em  ideal of the pair $(\nu, \mu)$}.
\end{definition}

The following powerful criterion was found by Kaliman and Kutzschebauch in \cite{KalimanKutzschebauch2008IM}.

%
%  CRITERION FOR COMPATIBLE PAIRS
%
\begin{theorem} \label{criterioncompatiblepair}
Let $X$ be an affine algebraic manifold on which the group of algebraic automorphisms 
$\Autalg(X)$ acts transitively. If there are compatible pairs $(\nu_i, \mu_i)$ and a point 
$p \in  X$ such that the vectors $\mu_i (p)$ form a generating set for $T_p X$, then $X$ 
has the algebraic density property.
\end{theorem}

Here we call a subset $ \{v_1, \ldots v_k\} \subset T_p X$ a {\em generating set} 
if the union of orbits of these vectors under the isotropy group $({\Autalg})_p (X)$ 
spans the tangent space $T_p X$.

The corresponding notion of a compatible pair, and the analogous theorem for general
(not necessarily algebraic) Stein manifolds was given by the same authors in 
\cite{KalimanKutzschebauch2015}. 

%
%   DEFINITION OF COMPATIBLE PAIRS ON STEIN MANIFOLDS
%
\begin{definition}\label{def:compatible}
A pair  $(\nu, \mu)$ of complete holomorphic vector fields on a Stein manifold $X$ is a 
{\em compatible pair} if the following conditions hold:
\begin{enumerate}[(1)] 
\item the closure of the linear span of the product of the kernels $\ker \nu \,\cdotp \ker \mu$ 
contains a nontrivial ideal $I \subset \Oscr(X)$, and
\item
there is a function $h \in \Oscr  (X)$ with $h\in \ker \mu$ and $\nu (h) \in \ker \nu\setminus \{0\}$. 
\end{enumerate}
The biggest ideal $I$ with this property is called the {\em ideal of the pair $(\nu, \mu)$}. 
If only condition (1) is satisfied, we call the pair semi-compatible.
\end{definition}

The following is an analogue of Theorem \ref{criterioncompatiblepair} for Stein manifolds.

\begin{theorem}
Let $X$ be a Stein manifold on which the group of holomorphic automorphisms $\Aut(X)$ act transitively. 
If there are compatible pairs $(\nu_i, \mu_i)$ such that there is a point $p \in  X$ where the vectors 
$\mu_i (p)$ form a generating set for $T_p X$, then $X$ has the density property.
\end{theorem}

This generalisation of the algebraic case % in Theorem \ref{criterioncompatiblepair} 
is crucial in proving that the Koras--Russell threefold has the density property, see (5) in the list below. 
Indeed, algebraic automorphisms do not act transitively on that threefold, 
whereas holomorphic automorphisms do act transitively.

Before we come to the complete list of examples of Stein manifolds known to have the density property, 
we would like to show the power of this criterion in some examples.

\begin{example}\label{Example:1}
On $\C^n$, $n\ge 2$, with coordinates $z=(z_1,z_2,\ldots, z_n)$ 
the pair of vector fields $(\frac {\partial} {\partial z_1}, \frac{\partial} {\partial z_2})$ 
is compatible, with the function $z_1$ and the ideal $I = \Oscr (\C^n)$ satisfying the conditions
in Definition \ref{def:compatible}. Since we can permute coordinates, 
$\{\frac{\partial} {\partial z_2}\}$ is a generating set for each tangent space. Thus, $\C^n$ has 
the density property.
\end{example}

\begin{example}\label{Example:2}
Denote an element of the special linear group $X= \SL_2 (\C )$ by 
\[
	A =  \left(\begin{array}{ccc}
	a_1 & a_2\\ b_1 & b_2 \end{array}\right). 
\]
The pair of vector fields on $X$ given by 
\[
	\delta_1 = b_1\frac{\partial}{\partial a_1}+b_2\frac{\partial}{\partial a_2}, \quad\ 
	\delta_2 = a_1\frac{\partial}{\partial b_1}+a_2\frac{\partial}{\partial b_2} 
\]
is compatible, with the ideal $I =\Oscr (X)$ and the function $h = a_1$.
Observe that the time-$t$ map of the field $\delta_1$ is adding $t$-times the first row to 
the second row, and vice versa for the field $\delta_2$. Hence these fields are 
tangent to $\SL_2 (\C)$. Since the adjoint representation of $\SL_2 (\C)$ is irreducible,
$\delta_2$ is a generating set at the identity. Thus, $\SL_2 (\C)$ has density property. 
The proof of the density property for $\SL_n (\C)$ and $\GL_n (\C)$,  $n\ge 3$, goes the same way.
\end{example}

\smallskip\noindent
\textbf{A list of examples of Stein manifolds known to have the density property:}

\smallskip\noindent
(1)  A homogeneous space $X=G/H$, where $G$ is a linear complex algebraic group and $H$ 
is a closed algebraic subgroup such that $X$ is affine and 
the connected component of $X$ is different from $\C$ and from $(\C^*)^n, n\ge 1,$ 
has the (algebraic) density property. 

It is known that if the subgroup $H$ is reductive then the space $X=G/H$ is always affine, but 
there is no known group-theoretic characterisation of $G$ which would say when is $X$ affine. 

The above result has a long history and includes all previously known examples from works of 
Anders\'en--Lempert, Varolin, T\'oth--Varolin, Kaliman--Kutzschebauch, and Donzelli--Dvorsky--Kaliman. 
The final result was obtained  by Kaliman and Kutzschebauch in \cite{KalimanKutzschebauch2017MA}. 
The manifolds $\C$ and $\C^*$ clearly do not have density property; however, the following 
problem is well known and seems notoriously difficult.

\begin{problem}\label{prob:Cstarn}
Does $(\C^*)^n$ for $n\ge 2$ have the density property? 
\end{problem}

It is conjectured that the answer is negative. More precisely, one expects that
all holomorphic automorphisms of $(\C^*)^n, n\ge 2,$ respect the volume form 
$\wedge_{i=1}^n \frac{dz_i}{z_i}$ up to a sign.

\smallskip\noindent
(2) The manifolds  $X$ given as a submanifold of $\C^{n+2}$ with coordinates $u\in \C$, $v\in \C$, $z\in \C^n$  by the equation $uv = p(z)$, where the zero fibre of the polynomial $p\in \C[\C^n]$
is smooth and reduced (otherwise $X$ is not smooth), have (algebraic)  density property; see
\cite{KalimanKutzschebauch2008MZ}. 

\smallskip
Before formulating the next result, recall that {\em Gizatullin surfaces} are by definition the normal affine
surfaces on which the algebraic automorphism group acts with an open orbit whose complement 
is a finite set of points. By the classical result of Gizatullin \cite{Gizatullin1971}, 
they are characterized by admitting a completion 
with a simple normal crossing chain of rational curves at infinity. 
Every Gizatullin surface admits a $\C$-fibration with at most one singular fibre which however
is not always reduced. Since the only affine algebraic 2-manifolds admitting semi-compatible pairs 
are $\C\times\C$ and $\C\times \C^*$ \cite{KalimanKutzschebauch2010IM}, 
the criterion from Theorem \ref{criterioncompatiblepair} is not applicable for surfaces, 
which makes the proof of the following result more cumbersome.

\smallskip\noindent
(3) Smooth Gizatullin surfaces which admit a $\C$-fibration with at most one singular and reduced fibre
have the density property (Andrist \cite{Andrist2018}). 
These surfaces are also called {\em generalised Danielewski surfaces}.
Special cases of this result were proved before in 
\cite{Donzelli,KalimanKutzschebauch2008MZ,AndristKutzschebauchPoloni2017}.

\smallskip\noindent
(4) The only known non-algebraic examples of Stein manifolds with the density property are, 
firstly, the holomorphic analogues of (2), namely, complex submanifolds $X$ of $\C^{n+2}$ 
with coordinates $u\in \C$, $v\in \C$, $z\in \C^n$ given
by an equation $uv = f(z)$, where the zero fibre of the holomorphic function $f \in \Oscr(\C^{n})$ is 
smooth and reduced (otherwise $X$ is not smooth); see \cite{KalimanKutzschebauch2008MZ}. 
Secondly, in the special case of (3) when the Gizatullin surface can be completed by four rational curves, 
the Stein manifolds given by the same algebraic equations (but using holomorphic functions 
as in the above example $uv=f(z)$) have the density property \cite{AndristKutzschebauchPoloni2017}.

\smallskip\noindent
(5) Certain hypersurfaces in $\C^{n+3}$ with coordinates 
$z=(z_0, z_1, \ldots , z_n) \in \C^{n+1}, x\in \C, y \in \C$,  
given by the polynomial equation $x^2y= a(z) + x b(z)$ where $\deg_{z_0} a \le 2$, $\deg_{z_0} b \le 1$ 
and not both degrees are zero, have the density property. 
(The exact conditions on $a$ and $b$ ensuring transitivity of the holomorphic automorphism group 
are rather technical.) This family includes the Koras--Russell threefold given in $\C^4$ by the equation 
\begin{equation}\label{eq:KR}	
	x +x^2y +s^2 +t^3 =0.
\end{equation}
This result of Leuenberger \cite{Leuenberger2016} is interesting in connection with the recognition problem for affine spaces; see Section \ref{sec:recognition}.

\smallskip\noindent
(6) The Calogero--Moser spaces $\mathcal{C}_n$, $n\in \N$, have the algebraic density property 
according to a recent result of Andrist \cite{Andrist2021}. The space $\mathcal{C}_n$ 
describes the phase space of a Calogero--Moser system, an $n$-particle system in classical physics 
with a certain Hamiltonian. %For the precise definition we refer the interested reader to the cited paper. 

%
%
%	VOLUME DENSITY PROPERTY
%   
%
\subsection{Volume density property}\label{sec:volumedensity}
The volume density property was the very first known density property (before the terminology was
introduced), which was discovered on Euclidean spaces $\C^n$ by Anders\'en \cite{Andersen1990}
in 1990. (This was the first paper on the subject and a precursor to Anders\'en's joint
work with Lempert \cite{AndersenLempert1992}.)  
We consider a holomorphic volume form $\omega$ on $X$, i.e., a nowhere vanishing  
holomorphic section of the canonical bundle 
$K_X=\wedge^n T^*X$, $n=\dim X$. (Note that $K_X$, being a line bundle over $X$,
must be topologically trivial for such a form to exist, and on a Stein manifold this necessary condition is 
also sufficient by the Oka--Grauert principle \cite[Theorem 5.3.1]{Forstneric2017E}.)
The Lie algebra $\ggot$ from Definition \ref{def:densitygeneral} is now the 
Lie algebra $\VF_\omega (X)$ of volume preserving holomorphic vector fields $\theta$, i.e., 
such that the Lie derivative of $\omega$ along $\theta$ vanishes, $L_\theta (\omega) = 0$. 
Let $\iota_\xi \omega$ denote the interior product of $\omega$ and a vector field $\xi \in \VFH (X)$.
Since $d\omega=0$, Cartan's formula for the Lie derivative gives
\begin{equation}\label{eq:Cartan}
	L_\theta (\omega)=d(\iota_\theta \omega) = \div_\omega(\theta) \omega.
\end{equation}
The function $\div_\omega(\theta)$ 
is called the divergence of $\theta$ with respect to $\omega$. Hence, $\VF_\omega (X)$ is the algebra of 
holomorphic vector fields on $X$ with vanishing divergence, $\div_\omega(\theta)=0$.

Since volume preserving vector fields do not form an $\Oscr$-module, the proof that the algebraic 
volume density property implies the holomorphic volume density property (see Kaliman and 
Kutzschebauch \cite{KalimanKutzschebauch2010IM}) is not straightforward.
For the same reason, the version of Theorem \ref{th:FR} of Forstneri\v c and Rosay \cite{ForstnericRosay1993}
for holomorphic automorphisms of $\C^n$ with Jacobian one (preserving the standard holomorphic 
volume form $dz_1\wedge \cdots\wedge dz_n$) requires the additional cohomological assumption 
$H^{n-1}(\Omega,\C)=0$ on the domain $\Omega$ of the isotopy.

The list of Stein manifolds known to enjoy the (algebraic) volume density property was rather short before 
an efficient criterion was established by Kaliman and Kutzschebauch in \cite{KalimanKutzschebauch2016TG}. 
We only state the holomorphic version from \cite{KalimanKutzschebauch2015}; the algebraic version is 
similar as in Theorem \ref{criterioncompatiblepair}. %We begin with some notation.
Let $X$ be a Stein manifold of dimension $n$ with a holomorphic volume form $\omega$. 
Denote by $\mathcal{Z}_{n-1}(X)$ the space of closed holomorphic differential $(n-1)$-forms on $X$.
The formula \eqref{eq:Cartan} shows that the map
\begin{equation}\label{eq:Theta}
	\Theta : \VF_\omega (X) \stackrel{\cong}{\to} \mathcal{Z}_{n-1}(X),
	\quad \ \xi \to  \iota_\xi \omega
\end{equation} 
is an isomorphism. By $\LieO (X)$ we denote the Lie subalgebra of $\VF_\omega (X)$ generated 
by complete holomorphic $\omega$-volume preserving vector fields on $X$.

%
%   KK CRITERION FOR STEIN MANIFOLDS WITH VOLUME DENSITY PROPERTY
%
\begin{theorem}  Let  $X$ be a Stein manifold with a holomorphic  volume form $\omega$.
Assume that there are pairs of divergence-free semi-compatible vector fields $(\xi_j, \eta_j)$ 
with ideals $I_j$ satisfying the following two conditions.
\begin{enumerate}[\rm (A)]
\item For every $x \in X$ the set $\{I_j (x) \xi_j (x) \wedge \eta_j (x) \}_{j=1}^k$
generates $\wedge^2 T_{x}X$.
\item The image of  $\Theta (\overline{\LieO (X)}) \subset \mathcal{Z}_{n-1} (X) $ under the
de Rham homomorphism $\Phi_{n-1} : \mathcal{Z}_{n-1} (X) \to H^{n-1} (X, \C)$ equals %coincides with 
$H^{n-1}(X , \C)$.
\end{enumerate}
Then $\Theta (\overline{\LieO (X)}) = \mathcal{Z}_{n-1} (X)$ and therefore
$\overline{\LieO (X)} = \VF_\omega (X)$, i.e., the manifold $(X,\omega)$ has the volume density property.
\end{theorem}

A highly nontrivial fact proved by Kaliman and Kutzschebauch 
\cite[Theorem 8]{KalimanKutzschebauch2015} using this criterion is that the product 
of two Stein manifolds $(X_1,\omega_1)$ and $(X_2,\omega_2)$ with the volume density property again 
has the volume density property for $\omega =\omega_1\wedge \omega_2$.

%
%   EXAMPLES OF STEIN MANIFOLDS WITH THE VOLUME DENSITY PROPERTY
%
\smallskip\noindent
\textbf{The list of examples of Stein manifolds known to have the volume density property:}

%The list of Stein manifolds for which the volume density property has been established up to now is as follows:

\smallskip\noindent
(1)  A homogeneous space $X=G/H$, where $G$ is a linear algebraic group and $H$ is a closed algebraic 
subgroup  such that $X$ is affine and admits a $G$-invariant volume form
(which is necessarily algebraic) has the (algebraic) volume density property
(Kaliman and Kutzschebauch \cite{KalimanKutzschebauch2017MA}). 
This includes earlier results of Varolin \cite{Varolin2001} and Anders\'en \cite{Andersen2000}
for the tori $(\C^*)^n$.

\smallskip\noindent
(2) The affine submanifolds $X$ in $\C^{n+2}$ with coordinates $u\in \C$, $v\in \C$, 
$z\in \C^n$, given by the equation $uv = p(z)$ where $p\in  \C[\C^n]$ is a polynomial whose 
zero fibre is smooth and reduced (to ensure that $X$ is not smooth), 
have (algebraic)  volume density property with respect to the unique algebraic volume form on them
(Kaliman and Kutzschebauch \cite{KalimanKutzschebauch2010IM}). Uniqueness 
%of the algebraic volume form 
(up to a multiplicative constant) follows from simple connectedness of the manifolds.

\smallskip\noindent
(3) The first, and up to now the only known non-algebraic examples
are certain holomorphic analogues of (2). Namely, for a nonconstant holomorphic function 
$f$ on $\C^n, n\ge 1,$ with $X_0=f^{-1} (0)$ reduced and smooth, and
such that the reduced cohomology $\wt H^{n-2}(X_0)$ vanishes if $n\ge 2$, 
the hypersurface $X\subset \C^{n+2}$ defined by $uv=f(z_1,\ldots,z_n)$ 
has the volume density property with respect to a certain volume form
(Ramos--Peon \cite{Ramos}).

\smallskip\noindent
(4) Certain hypersurfaces $X$ in $\C^{n+3}$ with coordinates 
$z=(z_0, z_1, \ldots , z_n) \in \C^{n+1}, x\in \C, y \in \C$ given by 
the equation $x^2y= a(z) + x b(z)$, where $\deg_{z_0} a \le 2$, $\deg_{z_0} b \le 1$ and not both 
degrees are zero, have volume density property for the volume form 
$\omega = \frac{dx}{ x^2} \wedge d z_0 \wedge d z_1 \wedge \ldots \wedge d z_n$
(Leuenberger \cite{Leuenberger2016}). 
This family includes the Koras--Russell threefold \eqref{eq:KR}. Again, the conditions ensuring 
transitivity of the group $\Aut_\omega(X)$ are rather technical.

\smallskip\noindent 
(5) The smooth fibres of the Gromov--Vaserstein fibration given by equations 
$
	p_n(z_1, \ldots, z_n) = a, 
$
where the polynomials $p_n\in  \C[\C^n]$ are defined inductively by 
\[
	p_{n+1} = z_{n+1} p_n + p_{n-1}, \quad p_0 = 1,  \quad p_1(z_1) =z_1, 
\]
have the volume density property for the unique algebraic volume form on them (De Vito \cite{Vito}). 
Uniqueness of the algebraic volume form (up to a constant) again follows from simple 
connectedness of the manifolds.

%
%  RELATIVE DENSITY PROPERTIES
%
\subsection{Relative density properties}\label{ss:relativedensity}
The next results concern density properties for  Lie algebras of holomorphic vector fields 
vanishing on subvarieties $Y \subset X$. It makes sense to consider such properties on a
Stein space $X$ when the subvariety $Y$ contains the singular locus $\Xsing $ of $X$. 
By using the same idea as in the proof of Theorem \ref{th:FR},
this density property naturally leads to theorems on approximation of isotopies of injective 
holomorphic maps $\Omega\to X$ on Runge domains $\Omega\subset X$,  
fixing $Y \cap \Omega$, by holomorphic automorphisms of $X$ fixing $Y$.

The following definitions were introduced in \cite{KutzschebauchLeuenbergerLiendo2015}.
Let $X$ be an affine algebraic variety, $Y\subseteq X$ be an algebraic subvariety 
containing $\Xsing$, and let $I=I(Y)\subseteq \C[X]$ denote the ideal of $Y$. 
Let $\VFalg(X,Y)$ be the $\C[X]$-module of vector fields vanishing on $Y$:
\[	
	\VFalg(X,Y)=\{\partial \in \VFalg(X) : \partial(\C[X])\subseteq I\}. 
\]
Let $\Liealg(X,Y)$ denote the Lie algebra generated by complete vector
fields in $\VFalg(X,Y)$.

\begin{definition} \label{def:ADP} 
{\rm (Assumptions and notation as above.)} 
\begin{enumerate}[\rm (a)] 
\item $X$ has the {\em strong algebraic density property relative to $Y$} 
if $\VFalg(X,Y) = \Liealg(X,Y)$.
\item $X$ has the {\em algebraic density property relative to $Y$} 
if there exists an integer $\ell\geq 0$ such that $I^\ell\VFalg(X,Y) \subseteq \Liealg(X,Y)$.
\end{enumerate}
\end{definition}

Note that condition (b) with $\ell=0$ is equivalent to condition (a).
We say that $X$ has the (strong) algebraic density property if these conditions hold for $Y=\Xsing$.

Except for the fact that we consider not necessarily smooth varieties,
the strong algebraic density property ($\ell=0$) is a version of Varolin's Definition \ref{def:densitygeneral}
of the density property for the Lie subalgebra $\ggot$ of vector fields
vanishing on $Y$. On the other hand, for $\ell>0$ our property in (b) is slightly weaker
than Varolin's definition since we generate the Lie subalgebra of
vector fields vanishing on $Y$ of order at least $\ell$ using complete
vector fields vanishing on $Y$ of possibly lower order than $\ell$. Still, 
this version of the algebraic density property  has the same remarkable
consequences as in Varolin's version of the algebraic density property for the group of 
holomorphic automorphisms of $X$ fixing the subvariety $Y$. 
In particular, the following analogue of Theorem \ref{th:FR} with interpolation on a subvariety holds.

%
%   APPROXIMATION WITH INTERPOLATION
%
\begin{theorem}\label{th:interpolationY}
Let $X$ be an affine algebraic variety and $Y\subset X$ be a closed algebraic subvariety
containing $\Xsing$. Assume that $\Omega\subset X$ is a Stein Runge domain
and $\Phi_t:\Omega\to X$ $(t\in [0,1])$ is an isotopy of injective holomorphic maps as 
in Theorem \ref{th:FR}. If condition (b) in Definition \ref{def:ADP} holds for an integer
$\ell\ge 0$ and the map $\Phi_t$ agrees with the identity  to order $\ell$ on $Y\cap \Omega$
for every $t\in[0,1]$, then $\Phi_1$ can be approximated uniformly on compacts in
$\Omega$ by holomorphic automorphisms of $X$ fixing $Y$ pointwise.
\end{theorem}

This result is often used to move a compact $\Oscr(X)$-convex subset of $X\setminus Y$
around by automorphisms of $X$ fixing $Y$. When $X=\C^n$ with $n\ge 2$,
this method can be used to find Fatou--Bieberbach domains containing $Y$ and omitting a certain 
compact set. (See Section \ref{sec:FB}.) Alternatively, approximating a contraction on a closed ball
$B\subset \C^n\setminus Y$ and fixing $Y$ gives a Fatou--Bieberbach domain in $\C^n\setminus Y$
containing $B$, thereby generalising Corollary \ref{cor:KL}. 
(See \cite[Corollary 4.12.2]{Forstneric2017E} and Theorem \ref{th:avoidvariety} in Section \ref{sec:FB}).  
In fact, this holds whenever $Y$ is an algebraic subvariety of $\C^n$ of codimension at least two
or, more generally, a tame subvariety of $\C^n$ 
of codimension at least two; see \cite[Definition 4.11.3 and Theorem 4.12.1]{Forstneric2017E}.

\begin{comment}
%
%   KALIMAN'S MESSAGE
%

The following result is due to Kaliman \cite[Theorem 2.17]{Kaliman2020} (2020).

\begin{theorem} 
Let $X$ be a complex affine flexible manifold and $Y$ be a closed complex subvariety of $X$ of  
codimension at least 2. Suppose that $X$ admits a pair of compatible vector fields.
Let $\Phi_t:\Omega\to \Phi_t(\Omega)\subset X\setminus Y$ $(t\in[0,1])$ 
be an isotopy of injective holomorphic maps between Stein Runge domains
with $\Phi_0=\Id$. Then $\Phi_1$ can be approximated uniformly on compacts in $\Omega$
by holomorphic automorphisms of $X$ which fix $Y$ pointwise.
\end{theorem}
\end{comment}

Here are the main examples when the relative  density property is known to hold.

\smallskip \noindent (1) 
The Euclidean space $\C^n$ for $n\ge 2$ has the relative algebraic density property 
with respect to any algebraic subvariety $Y\subset \C^n$ of codimension at least 2.
If moreover the dimension of the Zariski tangent space $T_y Y$ 
at every point of $Y$ is at most $n-1$, then $\C^n$ has the strong
algebraic density property with respect to $Y$ \cite[Theorems 4 and 6]{KalimanKutzschebauch2008IM};
this holds in particular if $Y$ is without singularities.
Hence, Theorem \ref{th:interpolationY} holds for any such subvariety $Y\subset \C^n$. 

% Franc: Is the same true for tame (=algebraically degenerate) subvarieties of codim at least two in C^n?

\smallskip\noindent (2)  Recall that a {\em locally nilpotent derivation} 
on an affine algebraic manifold is a complete algebraic vector field which generates an algebraic 
subgroup (isomorphic to $\C$ with $+$) of the algebraic automorphism group $\Autalg (X)$. 
An affine algebraic manifold $X$ is flexible in the sense of Arzhantsev et al.\  
\cite{ArzhantsevFlennerKalimanKutzschebauchZaidenberg2013DMJ} 
if locally nilpotent derivations span the tangent space  at every point of $X$.
Kaliman proved in \cite[Theorem 2.15 and Remark 2.16]{Kaliman2020} that
if $X$ is a flexible affine algebraic manifold with a compatible pair of locally nilpotent derivations, 
then $X$ has the density property relative to any algebraic subvariety $Y$ of codimension at least $2$. 
Again, it follows that Theorem \ref{th:interpolationY} holds for any such subvariety 
$Y\subset X$ (cf.\ \cite[Theorem 2.17]{Kaliman2020}). %By Example \ref{Example:1} 
This vastly generalizes the first part of (1), while Example \ref{Example:2} shows that 
$\SL_n (\C)$ has the density property relative to any algebraic subvariety of codimension 
at least $2$.

\smallskip \noindent (3) If $X$ is a normal affine toric variety of dimension $n\ge 2$ and 
$Y$ is a $\TT \cong (\C^*)^n$-invariant closed subvariety of $X$ containing
$\Xsing$, then $X$ has the algebraic density property relative to $Y$ if and only if
$X\setminus Y\neq \TT$ \cite[Theorem 3.7]{KutzschebauchLeuenbergerLiendo2015}. 
Affine toric surfaces with the strong algebraic density property were classified
by Kutzschebauch, Leuenberger, and Liendo \cite[Corollary 5.5]{KutzschebauchLeuenbergerLiendo2015}.
In the same paper, they gave an upper bound for the vanishing order $\ell$ in Definition \ref{def:ADP}.
This includes the results of Varolin  \cite{Varolin2001} for $Y$ a linear subspace of $X=\C^n$.
\subsection{Fibred density properties}\label{ss:fibreddensity}
Let $\pi : X \to B$ be a holomorphic map between complex manifolds. 
We say that $X$ has the {\em fibred density property 
with respect to $\pi$} if the Lie algebra $\mathfrak g$ of holomorphic vector fields $ \theta$ tangent 
to the fibres of $\pi$, i.e. fulfilling the condition $\mathrm{d}\pi (\theta) = 0$, has the density property; 
see Definition \ref{def:densitygeneral}.

The simplest case is a trivial fibration $W\times X \to W, \ (w,x) \to w$, where $X$ is a Stein manifold 
with the density property and the parameter space $W$ is a Stein manifold. For this case the fibred 
density property is easy to prove \cite{Kutzschebauch2005,KutzschebauchRamosPeon2017}.
The implication for the corresponding group of holomorphic automorphisms (in a simple situation
of a projection $\C^k\times \C^m\to\C^k$) is given by Theorem \ref{th:parameters} from the introduction.

The next case where the fibred density property is known is much more complicated.  
To formulate it, we need to recall the classical invariant-theoretic quotient of the action of 
$\mathrm{SL}_n (\C)$ on its Lie algebra $\mathfrak{sl_n}  (\C)$ of complex $n\times n$ matrices 
by conjugation (the adjoint representation). 
Denote by $\sigma_1, \dots, \sigma_n : \C^n \to \C$ the elementary symmetric polynomials in $n$ 
complex variables. Let $ \mathrm{Eig} : \mat{n}{n}{\C} \to \C^n$ assign to each matrix a vector of its 
eigenvalues. Denote by $\pi_1 := \sigma_1 \circ \mathrm{Eig}, \dots, \pi_n := \sigma_n \circ \mathrm{Eig}$
the elementary symmetric polynomials in the eigenvalues. By symmetrizing we avoid ambiguities 
caused by the order of eigenvalues in the definition of $\mathrm{Eig}$ and obtain a 
polynomial map $(\pi_1, \dots, \pi_n)$ such that 
$
	\chi_A(\lambda) = \lambda^n + \sum_{j=1}^n (-1)^j \pi_j(A)  \lambda^{n-j}
$
is the characteristic polynomial of the matrix $A$.

Since $\mathrm{trace} A = 0$ for  $A \in \mathfrak{sl_n}$, the map $\pi_1$, the sum of the eigenvalues, 
is the zero map.

Consider the fibration $\pi := (\pi_2, \dots, \pi_n) : \mathfrak{sl}_n \to \C^{n-1}$.  
A generic fibre, i.e.\ a fibre above a base point with no multiple eigenvalues, consists of exactly one 
equivalence class of similar matrices, so it is a homogeneous space of $\SL_n (\C)$ and hence smooth. 
A fibre above a base point with multiple eigenvalues decomposes into several strata of $\SL_n(\C)$ 
orbits, the largest one being the orbit of a matrix with the largest possible Jordan blocks. The %complicated 
structure of these fibres is well-studied in classical invariant theory, see e.g.\ Kraft \cite{Kraft1984}. 
The fibration $\pi$ is the invariant theoretic quotient $\mathfrak{sl}_n (\C) /\!\!/ \SL_n(\C)$, i.e., 
the algebra of $\SL_n(\C)$-invariant polynomial/holomorphic functions on $\mathfrak{sl}_n (\C)$ is 
exactly the pull-back of the algebra $\pi^* (\mathcal{ O} (\C^{n-1}))$ of holomorphic functions on $\C^{n-1}$.
The proof that $\pi$ has relative density property
(Andrist and Kutzschebauch \cite{AndristKutzschebauch2015P})  
uses a special family of complete holomorphic vector fields in this algebra, 
derived from certain one-parameter subgroups of $\SL_n (\C)$ exactly as the 
shears \eqref{eq:shear} on $\C^n$ are derived from the shear vector fields $\frac{\partial}{\partial z_i}$. 
In the same paper, the authors applied this density property to provide a set of generators for a dense 
subgroup of the automorphism group of the spectral ball.

%
%
%   SECTION: SYMPLECTIC DENSITY PROPERTY
%
\subsection{Symplectic density property}\label{ss:symplectic}
Let $X$ be a complex manifold of dimension $2n$ with a {\em holomorphic symplectic form} $\omega$, 
a closed holomorphic $2$-form whose highest exterior power $\omega^n$ is nowhere vanishing.
The standard holomorphic symplectic form on $\C^{2n}$ with coordinates 
$z_1, \ldots , z_n, w_1, \ldots, w_n$ is
\begin{equation}\label{eq:symplecticform}
	\omega = \sum_{i=1}^n dz_i \wedge dw_i.
\end{equation}
A symplectic form on a complex surface is the same thing as a holomorphic volume form.

Let $\Omega$ be a domain in a complex symplectic manifold $(X,\omega)$. A holomorphic map 
$f:\Omega\to X$ is called symplectic if $f^*\omega=\omega$. 
Assume that $\Phi_t:\Omega\to X$ is a smooth isotopy of injective symplectic holomorphic 
maps with the infinitesimal generator $V_t\in \VF(\Omega_t)$, $t\in [0,1]$.
Differentiating the identity $\omega= \Phi_t^*(\omega)$ on $t$ and taking into account
the Cartan formula $L_V\omega=d(\iota_V\omega) + \iota_V d\omega =d(\iota_V\omega)$ 
for the Lie derivative gives
\[
	0 = \frac{d}{dt} \Phi_t^*(\omega) = \Phi_t^*(L_{V_t} \omega) =
	\Phi_t^*(d( \iota_{V_t} \omega)),
\]
which holds if and only if $d(\iota_{V_t} \omega) =0$ for all $t\in[0,1]$.
A holomorphic vector field $V$ satisfying $d(\iota_V \omega)=0$
is called {\em symplectic}. This shows that flows of symplectic 
holomorphic maps are generated by symplectic vector fields, and vice versa. 
A vector field $V$ is called {\em Hamiltonian} if $\iota_V \omega$ is an exact $1$-form, and 
a holomorphic function $H\in \Oscr(X)$ satisfying 
$
	dH=\iota_V \omega
$
is the Hamiltonian of $V$. Conversely, every holomorphic function $H$ on $X$ determines
a Hamiltonian vector field $V_H$ by the above equation. On $\C^{2n}$
with coordinates $(z,w)$ and the symplectic form \eqref{eq:symplecticform}, every 
symplectic vector field is Hamiltonian of the form
\[
	V_H = \sum_{i=1}^n \frac{\di H}{\di w_j}\frac{\di}{\di z_j} 
		- \frac{\di H}{\di z_j}\frac{\di}{\di w_j},\quad\ H\in\Oscr(\C^{2n}). 
\]

The alternating bilinear form on $\C^{2n}$ defined by 
\begin{equation} \label{eqAC5.1}
   	\wt \omega(u, v)=\sum_{j=1}^n u_j v_{n+j} - u_{n+j}v_j, \quad \ u,v\in \C^{2n}
\end{equation}
is the {\em standard linear symplectic form} on $\C^{2n}$.
The corresponding differential form is $\omega$ \eqref{eq:symplecticform}.

The following algebraic density property was proved by Forstneri\v c 
\cite[Proposition 5.2]{Forstneric1996MZ} in 1996. This was the first known density property 
following the original ones of Anders\'en and Lempert, and it predates the formal 
introduction of this notion.

\begin{proposition} 
Let $\omega$ be the symplectic form \eqref{eq:symplecticform} on $\C^{2n}$, and let $\wt \omega$ be given by 
\eqref{eqAC5.1}. The Lie algebra of polynomial Hamiltonian vector fields on $(\C^{2n},\omega)$
is generated by the complete Hamiltonian vector fields of the form
\begin{equation} \label{eq:Hamiltonian}
	V(x) = f\left( \wt\omega(x, v)\right)  \sum_{i=1}^{2n} v_j\frac{\di}{\di x_j}, 
	\quad\ x\in \C^{2n},\ v\in \C^{2n},\ f\in \C[\C^{2n}]. 
\end{equation} 
\end{proposition}

The polynomial vector field \eqref{eq:Hamiltonian} generates the flow 
\begin{equation}\label{eqAC5.5}
		\Phi_t(x)= x+ t f\left( \omega(x, v)\right) v, \quad\ t\in\C,\ x\in\C^{2n}
\end{equation}
consisting of symplectic polynomial shear automorphisms of $(\C^{2n},\omega)$. 
The proof of Theorem \ref{th:FR} gives the following result 
(see \cite[Proposition 2.3 and Remark]{Forstneric1995JGEA}
and \cite[Theorem 5.1]{Forstneric1996MZ}). 

\begin{theorem} \label{th:ACTheorem5.1}
Assume that $\Omega$ is a pseudoconvex domain in $\C^{2n}$, $n\in \N$, with $H^1(\Omega,\C)=0$.
If $\Phi_t:\Omega\to \C^{2n}$ $(t\in [0,1])$ is a $\Cscr^1$ isotopy of injective symplectic 
holomorphic maps (with respect to the symplectic form \eqref{eq:symplecticform}) 
such that $\Phi_0$ is the identity map on $\Omega$ and the domain $\Omega_t=\Phi_t(\Omega)$ 
is Runge in $\C^{2n}$ for every $t\in[0,1]$, then $\Phi_1$ can be approximated uniformly on compacts 
in $\Omega$ by compositions of symplectic shears \eqref{eqAC5.5}.
In particular, the group generated by symplectic shears \eqref{eqAC5.5} is dense in the
group $\Aut_\omega(\C^{2n})$ of symplectic holomorphic automorphisms of $(\C^{2n},\omega)$.  
\end{theorem}

The condition $H^1(\Omega,\C)=0$ and pseudoconvexity of $\Omega$ ensure that every symplectic holomorphic vector field on $\Omega_t$ is Hamiltonian, and hence by the Runge condition on $\Omega_t$
it can be approximated uniformly on compacts in $\Omega_t$ by polynomial Hamiltonian
vector fields on $\C^{2n}$. (Like in the second part of Theorem \ref{th:FR}, 
pseudoconvexity is used in order to know that de Rham cohomology can be computed 
by means of holomorphic differential forms.) 

Of all density properties presented in this section, this symplectic density property 
is the only one that has not been developed yet for more general Stein symplectic manifolds, 
and no effective criteria are known (except on Stein surfaces where a symplectic form
is just a holomorphic volume form). The  problem in mimicking 
the criteria of Kaliman and Kutzschebauch for density or volume density properties lies in finding a module 
structure to apply the theory of coherent analytic sheaves, in particular, Theorems A and B. 
A  complete symplectic vector field remains complete when multiplied with a function in its kernel; 
however, it stays symplectic only if it is multiplied by a function $f(H)$ of the Hamiltonian $H$ 
of the vector field. For example, the symplectic density property would be much more interesting for the 
Calogero--Moser spaces than the density property proved by Andrist in \cite{Andrist2021}.

Since the holomorphic cotangent bundle of a complex manifold carries a natural holomorphic 
symplectic structure, it seems reasonable to propose the following problem.

\begin{problem}
Does the cotangent bundle of a Stein manifold with the density property enjoy 
the symplectic density property?
\end{problem}

%
%
%   SECTION: AUTOMORPHISMS WITH GIVEN JETS
%
%
\section{Automorphisms with given jets}\label{sec:jets}

The Anders\'en--Lempert theory is about Stein manifolds with large automorphism groups.
One aspect of this phenomenon is demonstrated by the approximation theorems mentioned 
in the introduction and the interpolation theorems discussed in Section \ref{sec:DP}.
One can go further and ask which jets of locally biholomorphic maps at 
a closed discrete set of points are jets of an automorphism.
The following result for finitely many points 
was proved by Anders\'en and Lempert \cite[Proposition 6.3]{AndersenLempert1992} 
and Forstneri\v c (see \cite[Proposition 2.1]{Forstneric1999JGEA} 
and \cite[Proposition 4.15.3]{Forstneric2017E}).  
The analogous result on Stein manifolds with the
density property is due to Varolin \cite{Varolin2000}. For symplectic holomorphic 
automorphisms of $\C^{2n}$, see L{\o}w et al.\ \cite{LowPereiraPetersWold2016}. 

\begin{proposition}\label{prop:Proposition2.1}
Let $m,n,N\in \N$, with $n>1$. Assume that
\begin{enumerate}[\rm (a)] 
\item $K$ is a compact polynomially convex set in $\C^n$,
\item $\{a_j\}_{j=1}^s$ is a finite set of points in $K$,
\item $p$ and $q$ are  points in $\C^n\setminus K$, and
\item $P: \C^n\to \C^n$ is a polynomial of order $m$ with nondegenerate linear part and $P(0)=0$.
\end{enumerate}
Given $\epsilon>0$, there exists $\Phi\in \Aut(\C^n)$ satisfying the following conditions:
\begin{enumerate}[\rm (i)] 
\item $\Phi(p)=q$ and $\Phi(z)=q+P(z-p)+O(|z-p|^{m+1})$ as $z\to p$,
\item $\Phi(z)=z+O(|z-a_j|^N)$ as $z\to a_j$ for each $j=1,2,\ldots,s$, and
\item[\rm (iii)] $|\Phi(z)-z|+|\Phi^{-1}(z)-z|<\epsilon$ for each $z\in K$.
\end{enumerate}
If in addition we have $JP(z)=1+O(|z|^m)$ as $z\to 0$ then there exists a polynomial 
automorphism $\Phi$ with $J\Phi\equiv 1$ satisfying conditions {\rm (i)--(iii)}.
\end{proposition}

These results have proved very useful in the construction of holomorphic automorphisms 
with interesting dynamical properties.

An inductive application of Proposition \ref{prop:Proposition2.1} leads to the following Mittag-Leffler 
interpolation theorem for automorphisms of $\C^n$ (Buzzard and Forstneri\v c
\cite[Theorem 1.1]{BuzzardForstneric2000}). Recall that a discrete sequence $a_j$ without
repetition in $\C^n$ is said to be {\em tame} if there is an automorphism $\Phi\in\Aut(\C^n)$
such that $\Phi(a_j)=(j,0,\ldots,0)\in \C\times \{0\}^{n-1}$ for all $j\in \N$. This notion 
was introduced and studied by Rosay and Rudin in \cite{RosayRudin1988}. 
(See also \cite[Section 4.6]{Forstneric2017E}.)

%
%
%  MAIN RESULT
%
%
\begin{theorem} \label{th:BF2000}
Assume that $n>1$, $a_j$ and $b_j$ $(j\in\N)$ are tame discrete sequences in $\C^n$ 
without repetitions, and $P_j:\C^n \to \C^n$ is a polynomial map with nondegenerate 
linear part and $P_j(0)=0$ for each $j\in\N$. Then there exists an automorphism $F$ of $\C^n$
such that for every $j = 1,2,\ldots$ we have $F(a_j)=b_j$ and
\[ 
        F(z)=b_j+P_j(z-a_j)+O(|z-a_j|^{m_j+1}), \quad  z\to a_j.
\] 
If in addition every $JP_j(z)=1+O(|z|^{m_j})$ for every $j$ and the sequences $a_j$ and $b_j$ 
are very tame, then there exists $F \in \Aut_1(\C^n)$ with these properties.
\end{theorem}

Points in a tame sequence can be permuted by automorphisms, and hence we can speak of tame
(closed) discrete sets. Rosay and Rudin \cite{RosayRudin1988} gave several criteria
for tameness, and they constructed nontame and even rigid discrete sets $A$ in $\C^n$
for any $n>1$, i.e., such that no nontrivial automorphism of $\C^n$ fixes $A$. There also
exist discrete sets whose complements are volume hyperbolic, meaning in particular
that every entire map $\C^n\to\C^n\setminus A$ has rank $<n$ at all points. 
(See also \cite[Secs.\ 4.6--4.7]{Forstneric2017E}.) 
An interesting use of such sets is shown in Sect.\ \ref{sec:twisted}.

The notion of a tame sequence was generalised to Stein manifolds with the density property 
in a couple of distinct ways by Andrist and Ugolini \cite{AndristUgolini2019} and 
Winkelmann \cite{Winkelmann2019}.

The following parametric version of Proposition \ref{prop:Proposition2.1}, due to 
Ramos--Peon and Ugolini \cite{RamosPeonUgolini2019}, 
generalises earlier results of Kutzschebauch and Ramos--Peon 
\cite{KutzschebauchRamosPeon2017} and Ugolini \cite{Ugolini2017}. 
It concerns interpolation of jets by holomorphic automorphisms at finitely many points of a Stein manifold 
$X$ with the density property, where the jets and the interpolating automorphisms of $X$ 
depend holomorphically on a parameter in another Stein manifold $W$.

\begin{theorem}\label{th:RPU}
Let $W$ and $X$ be Stein manifolds and suppose that $X$ has the density property. 
Let $k\ge 0$ and $N\ge 1$ be integers, and let $x_1,\ldots,x_N$ be distinct points in $X$. 
Let $Y$ denote the space of $N$-tuples $\gamma=(\gamma_1,\cdots,\gamma_N)$ of 
$k$-jets at $x_1,\ldots,x_N$, respectively, with nondegenerate linear parts and 
distinct values at $x_1,\ldots,x_N$. Given a null-homotopic holomorphic map
$\gamma=(\gamma_1,\cdots,\gamma_N):W\to Y$,
there exists a null-homotopic holomorphic map $F:W\to \Aut(X)$
such that the  $k$-jet of $F(w) \in\Aut(X)$ at $x_i$ is $\gamma(w)_i$ 
for $i=1,\ldots,N$ and all $w\in W$. 
\end{theorem}

\newpage

%
%
%   SECTION: FATOU--BIEBERBACH DOMAINS
%
%
\section{Fatou--Bieberbach domains}\label{sec:FB}

A domain $\Omega\subsetneq\C^n$ which is biholomorphic to $\C^n$ is called a 
{\em  Fatou--Bieberbach domain}. A biholomorphic map $F:\C^n\to\Omega$ onto such a domain
(and its inverse $F^{-1}$) is called a {\em Fatou--Bieberbach map}.
Every injective polynomial map $\C^n\to\C^n$ is an automorphisms with a polynomial inverse 
(see Rudin \cite{Rudin1995}), so there are no algebraic Fatou--Bieberbach maps.

All early constructions of Fatou--Bieberbach domains relied upon the theory of normal forms 
of local biholomorphisms at an attractive or a repelling fixed point. 
The following result has a complex genesis as explained in the sequel.

\begin{theorem} \label{th:RR}
Let $F:U\to F(U)$ be a biholomorphism on an open neighbourhood $U\subset\C^n$
of the origin such that $F(0)=0$ and the eigenvalues $\lambda_i$ of the differential $dF_0$ satisfy 
\begin{equation}\label{eq:ordered}
	1>|\lambda_1|\ge |\lambda_2|\ge \cdots\ge |\lambda_n|>0.
\end{equation}
After shrinking $U$ around the origin,
there is a biholomorphism $\psi:U\to \psi(U)\subset \C^n$ with $\psi(0)=0$ 
such that $G=\psi\circ F\circ \psi^{-1}$ is a  polynomial automorphism of $\C^n$ of the form 
\begin{equation}\label{eq:G}
	G(z) = Az + (0,g_2(z),\ldots,g_n(z)),\quad\ z\in\C^n,
\end{equation}
where $A$ is a lower-triangular matrix with the eigenvalues $\lambda_i$,
and every component $g_j(z)$ is a polynomial in the variables $z_1,\ldots,z_{j-1}$
containing no constant or linear terms. 
\end{theorem}

In fact, choosing $k\in\N$ such that $|\lambda_1|^k<|\lambda_n|$, $G$ may be chosen 
a polynomial map of degree $k$ such that every monomial
$z_1^{m_1}z_2^{m_2}\cdots z_n^{m_n}$ $(m_1+m_2+\cdots +m_n\ge 2)$ 
in a component $g_j(z)$ of $G$ \eqref{eq:G} is {\em resonant}:  
$\lambda_j=\lambda_1^{m_1}\lambda_2^{m_2}\cdots\lambda_n^{m_n}$.
Since the eigenvalues satisfy \eqref{eq:ordered}, it follows that 
$m_j=\cdots=m_n=0$ and hence $G$ is lower-triangular.

The first complete proof of Theorem \ref{th:RR} was given by Rosay and Rudin 
\cite[Appendix]{RosayRudin1988}; see also \cite[Sect.\ 4.3]{Forstneric2017E}. 
The result was claimed by Reich \cite{Reich1969-1,Reich1969-2} in 1969,
and it was used by Dixon and Esterle \cite{DixonEsterle1986} in 1986.
However, a gap in the proof of convergence of the normalization maps in
\cite{Reich1969-2} (even on the formal level) 
was pointed out by Rosay and Rudin \cite[p.\ 49]{RosayRudin1988}.
In the special case when the matrix $A$ of the differential $dF_0$ is diagonalisable and there are 
no resonances between the eigenvalues, we can conjugate $F$ locally near $0$ 
to the linear map $z\mapsto Az$ with $A=\mathrm{diag}(\lambda_1,\ldots,\lambda_n)$.
In the special case when all eigenvalues agree, this was proved 
by Poincar\'e \cite{Poincare1890} in 1890. Poincar\'e indicated
that this leads to the existence of injective holomorphic maps $\C^n\to\C^n$ with 
non-dense image (since there are automorphisms with several attracting fixed points),
although he did not provide a specific example. 
The formal normal form in the general case was developed by Leau \cite{Leau1897}
in 1987, who also obtained analytic solution under an extra technical hypothesis.
Further work was done by Picard \cite{Picard1900,Picard1905} in the period 1900-1905. 
Much later, Sternberg \cite{Sternberg1957} (1957) solved the normalization problem 
for smooth local diffeomorphisms of $\R^n$ at an attracting fixed point. 
His remark (see bottom of page 816 in his paper), that \cite[proof of Theorem 2]{Sternberg1957} 
also applies in the real analytic case, does not seem supported by any details in his paper.

Suppose now that $X$ is a complex manifold of dimension $n$ 
and $F:X\to X$ is an injective holomorphic map with an attracting 
fixed point at $p\in X$. Denote by $F^k$ the $k$-iterate of $F$ for $k\in \N$. The domain
\begin{equation}\label{eq:basin}
	\Omega_{F,p} = \Big\{x\in X: \lim_{k\to +\infty} F^{k}(x) =p \Big\}
\end{equation}
is called the {\em basin} of $F$ at $p$. Note that $\Omega_{F,p}$ is the increasing union of preimages 
$(F^k)^{-1}(V)$ for $k\in \N$, where $V\subset F(\C^n)$ is any neighbourhood of $p$.
By taking $V$ to be connected, we see that $\Omega_{F,p}$ is connected.
By Theorem \ref{th:RR} there are a neighbourhood $U\subset X$ of $p$, with $F(U)\subset U$, 
and a local chart $\psi:U\to \psi(U)\subset \C^n$ with $\psi(p)=0$ such that $G=\psi\circ F\circ \psi^{-1}$
is a lower-triangular polynomial map of the form \eqref{eq:G}. Clearly, 
such $G$ is an automorphism of $\C^n$ with a globally attracting fixed point at $0$, so 
$\Omega_{G,0}=\C^n$. 
Then, $G^k=\psi\circ F^k\circ \psi^{-1}$ for each $k\in\N$, which is equivalent to
$\psi=G^{-k}\circ\psi \circ F^k$. As $k\to\infty$, this defines a biholomorphic map 
$\Psi:\Omega_{F,p}\to \Omega_{G,0}=\C^n$, and hence shows the following.

%
%  BASIN OF AN INJECTIVE HOLOMORPHIC MAP AT AN ATTRACTIVE POINT 
%
\begin{theorem} \label{th:basin}
If $X$ is a complex manifold and $F:X\to X$ is an injective holomorphic map 
with an attracting fixed point $p\in X$, then the basin \eqref{eq:basin} 
is biholomorphic to $\C^n$, $n=\dim X$. 
\end{theorem}

%As explained above, this result has a complex history going back more than a century, with the punchline given by Rosay and Rudin \cite{RosayRudin1988} in 1988. 

When $n>1$, it may happen that the basin of an automorphism $F\in\Aut(\C^n)$ is not all of $\C^n$.
In fact, $F$ may have several (even countably many) attracting fixed points, 
and their basins are pairwise disjoint Fatou--Bieberbach domains in $\C^n$.

Before proceeding, we mention that fixed points $p\in\C^n$ of automorphisms $F\in\Aut(\C^n)$ 
which are hyperbolic, in the sense that the eigenvalues $\lambda_i$ of $dF_p$ satisfy
$|\lambda_i|\ne 1$, have also been studied. In favorable cases, $F$ can again be linearized at $p$.
In general there may be infinitely many resonances, which makes the study of the normal
form much more involved. There is an outstanding conjecture of Bedford that 
for a hyperbolic fixed point $p$ of $F\in\Aut(\C^n)$, the stable and the unstable manifolds
are biholomorphic to Euclidean spaces of appropriate dimensions. For a survey of 
this topic see Abbondandolo et al.\ \cite{Abbondandolo2014}. 

Using the idea behind Theorem \ref{th:basin}, Fatou \cite{Fatou1922-2} constructed in 1922 
birational self-maps of $\C^2$ whose images are not dense in $\C^2$. 
It was Bieberbach \cite{Bieberbach1933} who in 1933 found the first
known example of an injective holomorphic map $\C^2\to\C^2$ with Jacobian one
and non-dense image $\Omega$ which is Runge in $\C^2$. Bieberbach's example is also described
by Stehl\'e in \cite{Stehle1972} (1972), who used it to find a properly embedded holomorphic 
disc in $\C^2$. (The point is that if a complex line $\Lambda \subset \C^2$ intersects
a Runge domain $\Omega$ but 
is not contained in it, then any connected component of $\Lambda\cap \Omega$
is Runge in $\Lambda$, hence biholomorphic to the disc $\D$.)
See also the monograph by Bochner and Martin \cite[Sect. III.1]{BochnerMartin1948} (1948).
In 1971, Kodaira \cite{Kodaira1971} gave an example of an injective holomorphic 
map $\C^2\to\C^2$ with constant Jacobian omitting a complex line of $\C^2$.
In 1983, Nishimura \cite{Nishimura1983} gave such an example 
$F:\C^2\to\C^2$ which omits a neighbourhood $U$ of a complex line. In 
\cite{Nishimura1984}, Nishimura investigated the shape of $U$ for some specific $F$,
and he proved that there is no injective holomorphic map 
from $\C^2$ into itself with constant Jacobian whose image omits the union 
$E$ of two complex lines in $\C^2$ and a neighbourhood of a point of $E$. % (Theorem 3). 
This partially answers the (still open) question whether there is an injective holomorphic map
from $\C^2$ into itself omitting two complex lines:

\begin{problem} 
Is there a Fatou--Bieberbach domain in $\C^*\times\C^*$?
\end{problem}

The following result concerning basins of polynomial automorphisms of $\C^2$ was proved by
Bedford and Smillie  \cite{BedfordSmillie1991IUMJ} in 1991. Note that such an automorphism
has constant Jacobian determinant which is smaller than one in absolute value. 

%
%   Polynomial basins
%
\begin{theorem}
A polynomial basin $\Omega\subset\C^2$ intersects each algebraic curve 
$V\subset\C^2$ in a nonempty set with compact closure $\overline \Omega\cap V$.
On the other hand, the closure $\overline \Omega$ does not contain any closed one dimensional
complex subvarieties of $\C^2$. 
\end{theorem}

%
%   THE PUSH-OUT METHOD
%
In 1986, Dixon and Esterle \cite{DixonEsterle1986} introduced a
more general method for constructing Fatou--Bieberbach domains. The underlying idea is that 
for certain pairs of pairwise disjoint compact sets $K,L\subset\C^n$ with
polynomially convex union $K\cup L$ one can find holomorphic automorphisms $\Phi\in\Aut(\C^n)$
which are close to the identity map on one of the sets, say $K$, and they push
the second set $L$ far away. If $K$ and $L$ are convex, this is easily 
achieved by a shear, a fact which was also explored by Rosay and Rudin in \cite{RosayRudin1988}. 
In the more general case when $K$ is polynomially convex
and $L$ is starshaped (or holomorphically contractible), the same can be done by 
applying Theorem \ref{th:FR} (we squeeze $L$ within itself almost to a point, slide it
far away from $K$, and approximate the final map by an automorphism).
An inductive application of this construction yields a sequence
$\Phi_j\in\Aut(\C^n)$ $(j\in\N)$ such that the sequence of compositions
$F_j=\Phi_j\circ\Phi_{j-1}\circ \cdots\circ\Phi_1$ converges on a neighbourhood of $K$ 
and diverges to infinity on $L$ as $j\to\infty$. If in addition the $\Phi_j$'s approximate 
the identity map sufficiently closely on an increasing sequence of compacts exhausting $\C^n$, 
then the domain of convergence $\Omega$ of the sequence $F_j$ is a Fatou--Bieberbach domain, and
the limit $F=\lim_{j\to\infty} F_j:\Omega\to\C^n$ is a biholomorphic map of $\Omega$
onto $\C^n$ such that 
\[ 
	K\subset \Omega\subset \C^n\setminus L.
\] 
A similar argument yields a Fatou--Bieberbach domain $\Omega'\subset \C^n$
with $L\subset \Omega'\subset \C^n\setminus K$. See 
Forstneri\v c and Ritter \cite[Proposition 9]{ForstnericRitter2014} 
or \cite[Proposition 4.4.4]{Forstneric2017E} for a precise statement.
Taking $L$ to be a finite set gives the following corollary.

\begin{corollary}\label{cor:KL}
Given a compact polynomially convex set $K\subset \C^n$, $n>1$, and a finite set 
$L\subset\C^n\setminus K$ there is a Fatou--Bieberbach domain $\Omega$
satisfying $K\subset \Omega\subset \C^n\setminus L$. %\eqref{eq:KL}.
\end{corollary}

This result was used by Forstneri\v c \cite{Forstneric2003AM} in his construction 
of holomorphic functions without critical points on any Stein manifold. More generally,
Theorem \ref{th:avoidvariety} below is used in his proof of the basic h-principle for holomorphic 
submersions $X\to\C^q$ from any Stein manifold with $\dim X>q\ge 1$, given in the same paper. 

This {\em push-out method} of Dixon and Esterle also became known as {\em random iteration}, 
since we are not iterating an automorphism but composing with a new one at every step. 
It was used even before the Anders\'en--Lempert theory, or without using it, 
to provide examples of Fatou--Bieberbach domains 
with interesting properties and to solve various problems.
We refer to the papers by Rosay and Rudin \cite{RosayRudin1988},
Forn{\ae}ss and Sibony \cite{FornaessSibony1992}, Globevnik and Stens{\o}nes
\cite{GlobevnikStensones1995}, Globevnik \cite{Globevnik1997,Globevnik1998}, 
Stens{\o}nes \cite{Stensones1997}, Forn{\ae}ss and Stens{\o}nes \cite{FornaessStensones2004}, 
among others. Globevnik and Stens{\o}nes \cite{GlobevnikStensones1995}
used random iterations of shears in coordinate directions to show
that every planar domain bounded by finitely many Jordan curves admits a proper holomorphic
embedding into $\C^2$. This was a major advance on the open problem asking which 
open Riemann surfaces admit a proper holomorphic embedding as
a closed complex curve in $\C^2$ (see Section \ref{sec:RS}).  
Random iterations of shears were also used by Stens{\o}nes \cite{Stensones1997}
(1997) in her construction of Fatou--Bieberbach domains having $\Cscr^\infty$ 
smooth boundaries. Note that a smoothly bounded Fatou--Bieberbach domain in $\C^2$
has Levi-flat boundary foliated by complex curves.
Related results of Globevnik \cite{Globevnik1997,Globevnik1998}
give Fatou--Bierbach domains with $\Cscr^1$ boundaries and with additional geometric
control on their location. Whether there exist Fatou--Bieberbach domains with 
real analytic boundaries remains an open problem. 

A result similar to the one of Globevnik \cite{Globevnik1998} 
was used by Buzzard and Hubbard \cite[Lemma 3.1]{BuzzardHubbard2000}
to show the following \cite[Theorem 4.1]{BuzzardHubbard2000}.

\begin{theorem}\label{th:BH}
For every algebraic subvariety $A$ of codimension at least two in $\C^n$, $n\ge 2$,
there exists a Fatou--Bieberbach domain $\Omega\subset\C^n$ such that 
$\overline \Omega\cap A=\varnothing$.  
\end{theorem}

With the exception of Corollary \ref{cor:KL}, the results mentioned so far 
were obtained by elementary constructions using shears. The Anders\'en--Lempert theory provides 
much more general construction methods. For example, we have the following result.

\begin{theorem}\label{th:avoidvariety}
Let $A$ be an algebraic subvariety of $\C^n$, $n\ge 2$, of codimension at least $2$.
Given a compact, polynomially convex and holomorphically contractible set 
$K\subset \C^n\setminus A$, there is a Fatou--Bieberbach domain $\Omega\subset \C^n$ 
with $K\subset \Omega$ and $\overline\Omega\cap A=\varnothing$.
\end{theorem}

Here is a sketch of proof. By the assumption there are a neighbourhood 
$U\subset \C^n\setminus A$ of $K$ and an isotopy of biholomorphic maps 
$\phi_t:U\to \C^n\setminus A$ $(t\in[0,1])$ such that $\phi_t(K)\subset K$ for all $t$
and $\phi_1(K)$ lies in an arbitrarily small closed ball $B$ around a point $p\in K$.
Clearly, the sets $K_t=\phi_t(K)$ are polynomially convex. 
Let $\Omega'\subset \C^n\setminus A$ be a Fatou--Bieberbach domain
such that $\overline{\Omega'}\cap A=\varnothing$ (see Theorem \ref{th:BH}). 
If $B$ is small enough, we can slide it into $\Omega'$ by an isotopy of translations
which keep the image of $B$ in $\C^n\setminus A$. Applying 
Theorem \ref{th:interpolationY} to the combined isotopy provides
an automorphism $\Phi$ of $\C^n$ fixing $A$ such 
that $\Phi(K)\subset \Omega'$. Then, $\Omega=\Phi^{-1}(\Omega')$
is a Fatou--Bieberbach domain satisfying the conclusion of Theorem \ref{th:avoidvariety}.

Starting from 2005, one of the main contributors of developments on 
Fatou--Bieber\-bach domains and their applications has been Wold
with collaborators. In his first paper \cite{Wold2005}, Wold showed the following results. 
\begin{enumerate}[\rm 1.]
\item 
For any $m\in\N\cup\{\infty\}$ there exist $m$ pairwise disjoint Fatou--Bieberbach domains
$\Omega_i$ such that any point $p\in \C^n\setminus \bigcup_{i=1}^m \Omega_i$
lies in the boundary of every $\Omega_i$.
\item
If $\{L_i\}_{i\in\N}$ is a collection of affine subspaces of $\C^n$ $(n>1)$, then there
exists a Fatou--Bieberbach domain $\Omega$ such that for every $i$, $\Omega\cap L_i$
is connected and $L_i\setminus \Omega \ne \varnothing$.
\item 
If $\{V_i\}_{i\in\N}$ is a collection of closed proper complex subvarieties of $\C^n$ 
$(n>1)$ then there exists a Fatou--Bieberbach domain $\Omega$ containing $\bigcup_{i} V_i$.
\end{enumerate}

An exciting development was the following result of Wold \cite{Wold2008} from 2008.

\begin{theorem}\label{th:nonRungeFB}
For any $n>1$ there exists a non-Runge Fatou--Bieberbach domain in $\C^n$.
\end{theorem}

Note that such domains cannot be limits of automorphisms of $\C^n$. 
The idea behind the construction is the following. Start with a Fatou--Bieberbach domain
$\Omega$ contained in $\C^*\times \C$. By Stolzenberg \cite{Stolzenberg1966MA} there is a compact set 
$K=M_1\cup M_2 \subset \C^*\times \C$, consisting of a pair of disjoint closed totally real discs 
$M_1$ and $M_2$, such that $K$ is $\Oscr(\C^*\times \C)$-convex but its polynomial hull
$\wh K$ contains the origin of $\C^2$. Since the Stein domain 
$\C^*\times \C$ has the density property, Theorem \ref{th:FR2} can be used to find a 
holomorphic automorphisms $\Phi$ of $\C^*\times \C$ such that 
$\Phi(K)\subset \Omega$. (It suffices to construct an isotopy which shrinks each of the discs $M_1,M_2$
almost to a point and then slide the new small discs into $\Omega$ within 
$\C^*\times \C$ such that the isotopy consists of $\Oscr(\C^*\times \C)$-convex sets.)
Therefore, $\Omega'=\Phi^{-1}(\Omega)\subset \C^*\times \C$ is a Fatou--Bieberbach domain
containing $K$ but not its polynomial hull $\wh K$. 
It follows that $\Omega'$ is not Runge in $\C^2$, although it is Runge in $\C^*\times\C$.

An interesting application of this result was Wold's construction \cite{Wold2010} in 2010 of the first known 
example of a non-Stein long $\C^2$. We describe these developments in Section \ref{sec:long}.

In 2010, Baader, Kutzschebauch and Wold \cite{BaaderKutzschebauchWold2010}
used Fatou--Bieberbach domains to construct the first known example of a knotted properly 
embedded holomorphic disc in $\C^2$. Their result was motivated by the problem, 
raised by Kirby, whether proper holomorphic embeddings of $\C$ or the unit disc into $\C^2$ 
can be topologically knotted. While the first problem remains open, the second one was solved in 
the affirmative in \cite{BaaderKutzschebauchWold2010}. The proof uses well-behaved 
Fatou--Bieberbach domains in $\C^2$ constructed by Globevnik in \cite{Globevnik1998},
containing small perturbations of the bidisc, and the existence of knotted holomorphic discs 
in the bidisc. It is unknown whether the disc admits an unknotted
proper holomorphic embedding into $\C^2$.

In 1998 Globevnik \cite{Globevnik1998} constructed Fatou--Bieberbach domains in $\C^n$ 
whose closures intersect the complex line $\C\times \{0\}^{n-1}$ in closed connected 
and simply connected domains which are arbitrarily small perturbations of the closed unit disc.
In 2012, Wold \cite{Wold2012} found a Fatou--Bieberbach domain in $\C^2$ whose intersection
with $\C\times \{0\}$ contains the unit disc as a connected component, thereby
answering a question of Rosay and Rudin \cite{RosayRudin1988}.

In 2015, Forstneri\v c and Wold \cite{ForstnericWold2015}
constructed Fatou--Bieberbach domains in $\C^n$, $n>1$, which contain a given compact set $K$ 
and avoid a totally real affine subspace $L\subset\C^n$ with $\dim_\R L<n$ 
such that $K\cup L$ is polynomially convex. 
This was used in \cite{ForstnericWold2015} to show that $\C^n\setminus L$ has 
certain Oka properties. Due to a recent result of Kusakabe it is now clear that such domains 
are Oka manifolds for $n\ge 4$, as follows easily from Theorem \ref{th:E}. 

A recent development is the following result of Forstneri\v c and Wold 
\cite[Theorem 1.1]{ForstnericWold2020MRL} from 2020. 
An interesting application is given by Theorem \ref{th:KOka}.  

%
%   MAIN THEOREM
%
\begin{theorem}\label{th:FW2020}
Let $K$ be a compact polynomially convex set in $\C^n$ for some $n>1$, 
$L$ be a compact polynomially convex set in $\C^N$ for some $N\in\N$, 
and $f:U\to \C^n$ be a holomorphic map
from an open neighbourhood $U\subset \C^N$ of $L$ such that $f(z)\in \C^n\setminus K$
for all $z\in L$. Then there are an open neighbourhood $V\subset U$ of $L$
and a holomorphic map $F:V\times \C^n\to \C^n$ such that for every $z\in V$,
$F(z,0)=f(z)$ and the map $F(z,\cdotp):\C^n\to \C^n\setminus K$ is injective.
\end{theorem}

It follows that $\Omega_z:=\{F(z,\zeta):\zeta\in\C^n\}$ is a Fatou-Bieberbach domain 
in $\C^n\setminus K$ with centre $F(z,0)=f(z)$ and depending holomorphically on $z\in V$.

The proof uses Anders\'en-Lempert theory with parameters.
It also applies to variable fibres $K_z\subset \C^n$ $(z\in L)$ with polynomially convex graph 
(see \cite[Remark 2.2]{ForstnericWold2020MRL}). 
For a convex parameter space $L\subset \C^N$ the analogous result 
holds if we replace $\C^n$ by an arbitrary Stein manifold having the density property
(see \cite[Theorem 3.1]{ForstnericWold2020MRL}). 

%
%  KALIMAN
%
So far we have focused on Fatou--Bierbach domain in Euclidean spaces $\C^n$, $n\ge 2$.
However, such domains abound in any Stein manifold with the density property
as was already observed by Varolin \cite{Varolin2000,Varolin2001}.
As an example, we mention the following recent result 
of Kaliman \cite[Corollary 2.18]{Kaliman2020} 
which generalizes Theorem \ref{th:avoidvariety}. 
Recall that an affine manifold $X$ is flexible in the sense of Arzhantsev et al.\ 
\cite{ArzhantsevFlennerKalimanKutzschebauchZaidenberg2013DMJ} 
if locally nilpotent derivations on $X$ span the tangent space at every point
(see Example (2) in Subsect.\ \ref{ss:relativedensity}). 

\begin{theorem} \label{th:Kaliman2.18}
Let $X$ be a complex affine flexible manifold and $Y$ be a closed algebraic subvariety of $X$ of  
codimension at least 2. Suppose that $X$ admits a pair of compatible vector fields
(see Definition \ref{def:compatible}). 
Then, every point $x\in X$ has a neighbourhood $\Omega\subset X\setminus Y$ 
which is biholomorphic to $\C^n$ with $n=\dim X$.
\end{theorem}

\newpage
%
%
%   SECTION: CONNECTIONS TO OKA THEORY
%
%
\section{Twisted complex lines in $\C^n$ and nonlinearizable automorphisms}\label{sec:twisted}

It is well known that a generic affine algebraic hypersurface $A\subset\C^n$ of sufficiently
large degree is (Kobayashi) hyperbolic and has hyperbolic complement $\C^n\setminus A$
(see Brotbek \cite{Brotbek2017}).  
Such $A$ is necessarily topologically complicated. On the other hand, proper polynomial embeddings
$F:\C^k\hra\C^{k+1}$ are believed to have non-hyperbolic complements $\C^{k+1}\setminus F(\C^k)$.
In particular, it was shown by Suzuki \cite{Suzuki1974} (1974) 
and Abhyankar and Moh \cite{AbhyankarMoh1975} (1975) 
that for every polynomial embedding $F:\C\hra \C^2$ there is a 
polynomial automorphism $\Phi$ of $\C^2$ such that $\Phi\circ F(\C)=\C\times \{0\}$.
It is therefore of  interest to know that there are proper holomorphic embeddings
whose complements are hyperbolic; see Theorem \ref{th:hyperboliccomplement}.

On the way to this result, we begin with the following application of Theorem \ref{th:FR} due to  
Forstneri\v c, Globevnik, and Rosay \cite{ForstnericGlobevnikRosay1996}.

\begin{theorem}\label{th:embedding-interpolation}
For every closed discrete set $B\subset \C^2$ there is a properly 
embedded complex line $F:\C\hra \C^2$ such that $B\subset F(\C)$.
\end{theorem}

The analogous result holds for embedding $\C\hra\C^n$, $n\ge 3$. This was proved  
in a more precise form (with interpolation on a pair of discrete sets) by Rosay and Rudin 
\cite[Theorem I]{RosayRudin1993} in 1993. However, their result for $n\ge 3$ is a very special case 
of the theorem, due to Acquistapace, Broglia, and Tognoli
\cite{AcquistapaceBrogliaTognoli1975} (1975), that for any Stein 
manifold $X$ of dimension $m\ge 1$, integer $n\ge 2m+1$, and proper
holomorphic embedding $\phi:X' \hra\C^n$ of a closed complex subvariety $X'$ of $X$
there is a proper holomorphic embedding $F:X\hra \C^n$ with $F|_{X'}=\phi$. 
(See also \cite[Theorem 9.5.5]{Forstneric2017E}.) The Rosay--Rudin theorem mentioned 
above amounts to the special case with $X'$ a closed discrete subset of $X=\C$.
However, for $n=2$ the methods in both mentioned papers only provide immersions $\C\to\C^2$
with the desired property (cf.\ \cite[Theorem II]{RosayRudin1993}). In this lowest dimensional case, 
and in the proof of the more general result in Theorem \ref{th:hyperboliccomplement} for $n\le 2k$, 
the use of Anders\'en--Lempert theory is essential.

The main idea behind the proof of Theorem \ref{th:embedding-interpolation} 
is to inductively twist a properly embedded complex line
$\C\hra \C^2$ such that it contains more and more points of the discrete set 
$B=\{b_1,b_2,\ldots\}$, and the sequence of embeddings converges to a proper holomorphic embedding. 
In the inductive step we are given a proper holomorphic embedding $F_k:\C\hra\C^2$ 
such that $F_k(\C)$ contains the first $k$ points $b_1,\ldots, b_k\in B$ but it does not contain the remaining
points of $B$ (the latter condition is easily arranged by a general position argument). 
Choose a disc $\Delta\subset\C$ with $\{b_1,\ldots, b_k\}\subset F_k(\Delta)$
and a compact ball $L\subset\C^2$ such that 
\[
	F_k(b\Delta) \cup \{b_{k+1},b_{k+2},\cdots\} \subset \C^2\setminus L.
\]
The union $K:=L\cup F_k(\overline\Delta)$ is then a compact polynomially convex set in $\C^2$.
Theorem \ref{th:FR} furnishes an automorphism $\Phi_k\in\Aut(\C^2)$ 
which is close to the identity map on $K$, it fixes the points
$b_1,\ldots, b_k$, and such that $b_{k+1}\in \Phi_k\circ F_k(\C)$. This gives the next 
embedding $F_{k+1}=\Phi_k\circ F_k:\C\hra\C^2$. 
An inductive application of this technique gives a sequence of embeddings such that 
$F=\lim_{k\to\infty}F_k:\C\hra\C^2$ is an embedding satisfying the conclusion of the theorem.
The same proof applies for any $n\ge 2$.

A related result of Buzzard and Forstneri\v c in \cite{BuzzardForstneric1997} (1997) 
yields Carleman approximation with interpolation on $\R\subset \C$ 
of proper holomorphic embeddings $\C\hra\C^n$ for any $n>1$.

Rosay and Rudin proved in \cite{RosayRudin1988} that for any $n>1$ there exist discrete sets 
$B\subset \C^n$ which cannot be mapped into an affine line by any holomorphic 
automorphism of $\C^n$; such sets are called non-tame. Applying 
Theorem \ref{th:embedding-interpolation} (and its generalisation to any $n>1$ mentioned above)
with such a set $B$ gives the following corollary.

%
%   NONSTRAIGHTENABLE COMPLEX LINES
%
\begin{corollary}\label{cor:nonstraightenable}
For any $n>1$ there is a properly embedded complex line $F:\C\hra \C^n$ 
such that no holomorphic automorphism of $\C^n$ maps $F(\C)$ onto an affine complex line. 
In other words, $F(\C)$ is not straightenable by automorphisms of $\C^n$.
\end{corollary}

This is in strong contrast with the result of Suzuki \cite{Suzuki1974} (1974) 
and Abhyankar and Moh \cite{AbhyankarMoh1975} (1975) 
that every polynomial holomorphic embedding $\C\hra \C^2$ is straightenable by a 
polynomial automorphism of $\C^2$. We refer to \cite[Sect.\ 4.18]{Forstneric2017E} for a 
survey of further results on this subject. This is an  example where the answer
to a certain problem in the holomorphic category differs from the answer in the 
algebro-geometric category.

A stronger result in the same spirit was obtained by Buzzard and Forn\ae ss 
\cite{BuzzardFornaess1996}, who constructed a proper holomorphic embedding $F:\C\hra \C^2$ such that 
$\C^2\setminus F(\C)$ is Kobayashi hyperbolic. Since the complement of an affine line in $\C^2$
is biholomorphic to $\C^*\times \C$ and hence is not hyperbolic, their result also implies 
Corollary \ref{cor:nonstraightenable}. To prove it, they constructed properly embedded complex lines
$\C\hra\C^2$ which contain arbitrarily small deformations of a well-chosen
discrete family of closed affine discs in $\C^2$ with Kobayashi hyperbolic complement.
If the approximations are close enough then the complement of the embedded line is
also hyperbolic. Their proof uses the same general idea as the proof of 
Theorem \ref{th:embedding-interpolation}, but the technical details are more involved. 
This line of results was developed further by Forstneri\v c \cite{Forstneric1999JGEA} 
and Borell and Kutzschebauch \cite{BorellKutzschebauch2006} who proved the following.

%
%  HYPERBOLIC COMPLEMENTS
%
\begin{theorem}\label{th:hyperboliccomplement}
For every pair of integers $1\le k<n$ there is a proper holomorphic embedding $F:\C^k\hra\C^n$ 
such that  $\C^n\setminus F(\C^k)$ is $(n-k)$-hyperbolic in the sense of Eisenman.
In particular, any entire map $\C^p\to \C^n\setminus F(\C^k)$ $(p\in\N)$ has rank less than
$n-k$ at each point. 
\end{theorem}

The theorem of Buzzard and Forn\ae ss \cite{BuzzardFornaess1995} corresponds to the case $k=1,\ n=2$.

As an application of Theorem \ref{th:embedding-interpolation}, 
Derksen and Kutzschebauch \cite{DerksenKutzschebauch1998} 
showed the following result which answered a long-standing open question.

%
%   NONLINEARIZABLE AUTOMORPHISM 
%
\begin{theorem} \label{th:nonlinearizable}
For every integer $n\ge 2$ there exists a nonlinearizable periodic holomorphic automorphism 
of period $n$ on $\C^{2+n}$. In particular, there is a nonlinearizable holomorphic involution on $\C^4$.
\end{theorem}

An outline of their proof can also be found in \cite[Sect.\ 4.19]{Forstneric2017E}.
The problem regarding the existence of nonlinearizable periodic automorphisms 
remains open on $\C^2$ and $\C^3$. A recent survey of the linearization problem 
for holomorphic automorphisms is available in \cite{Kutzschebauch2020}.

%
%
%   SECTION: EMBEDDING OPEN RIEMANN SURFACES IN C^2
%
%
\section{Embedding open Riemann surfaces in $\C^2$}\label{sec:RS}

It has been known since mid-1950s that every Stein manifold $X$ of dimension $n$  
embeds as a closed complex submanifold in a Euclidean space $\C^{2n+1}$. 
We refer to \cite[Sect.\ 2.4]{Forstneric2017E} or \cite[Sect.\ 2]{Forstneric2018Korea} 
for a discussion and references to the early works on the subject. 
The smallest possible value of $N$ for any $n>1$ was found by Eliashberg and Gromov 
\cite{EliashbergGromov1992AM} and Sch\"urmann \cite{Schurmann1997}
who showed that every Stein manifold $X$ of dimension $n\ge 1$ immerses properly holomorphically in
$\C^M$ with $M = \left[\frac{3n+1}{2}\right]$, and if $n>1$ then $X$ embeds properly holomorphically
in $\C^N$ with $N = \left[\frac{3n}{2}\right]+ 1$. Their proof relies on an application of the Oka principle 
for sections of holomorphic fibre bundles with Oka fibres over a Stein manifold. 
A complete exposition can also be found in \cite[Sections 9.3--9.4]{Forstneric2017E}.

The proof of this embedding theorem breaks down in the lowest dimensional case $n=1$
(i.e., $X$ is an open Riemann surface) and $N=2$. The following 
{\em Forster--Bell--Narasimhan Conjecture} \cite{Forster1970,BellNarasimhan1990}
is one of the oldest open problems in complex analysis.

\begin{problem}\label{prob:RS}                                  
Does every open Riemann surface embed properly holomorphically in $\C^2$?
\end{problem}

A history of the rather sporadic progress on this problem can be found in 
\cite[Sect.\ 9.10]{Forstneric2017E}. We mention in particular that 
Globevnik and Stens{\o}nes \cite{GlobevnikStensones1995} proved in 1995 
that every finitely connected domain in $\C$ without isolated boundary points embeds properly 
holomorphically into $\C^2$. Their proof uses Fatou--Bieberbach domains, constructed 
as domains of convergence of random sequences of shears in coordinate directions. 
Further progress using the same techniques was made by 
\v Cerne and Globevnik \cite{CerneGlobevnik2000}
and \v Cerne and Forstneri\v c \cite{CerneForstneric2002}.

A new method based on Anders\'en--Lempert theory was introduced into the subject
a decade later by Wold \cite{Wold2006MZ,Wold2006IJM,Wold2007}. 
Assume that $M$ is a compact bordered Riemann surface
(every such is conformally equivalent to a domain in a compact Riemann surface obtained 
by removing finitely many pairwise disjoint discs \cite[Theorem 8.1]{Stout1965TAMS})
and $F:\overline M\hra \C^2$ is a smooth embedding which is holomorphic on $M$. 
We wish to show that this embedding can be modified so that the boundary curves diverge to infinity 
while the interior $M\setminus bM$ becomes embedded in $\C^2$ as a closed complex curve.
The main idea introduced by Wold is the following. 
Write $bM=\bigcup_{i=1}^m C_i$ where each $C_i$ is a smooth closed curve.
Assume in addition that each $C_i$ contains a point $p_i$ such that 
the affine complex line in $\C^2$ through the point $F(p_i)=(a_i,b_i)\in\C^2$ in 
the second coordinate direction intersects $F(M)$ only at $F(p_i)$. 
Such point $p_i$ is said to be {\em exposed} by the map $F$. We apply to $F$
a rational shear of the form
\[
	G(z_1,z_2)=\Big(z_1,z_2+ \sum_{i=1}^m \frac{c_i}{z_1-a_i} \Big)
\]
for a suitable choice of the numbers $c_i\in\C^*$. Each point $F(p_i)$ 
is sent to infinity, $G$ has no other poles on $F(M)$,
and the surface $\Sigma=(G\circ F) (M\setminus \{p_1,\ldots,p_m\}) \subset \C^2$ 
is holomorphically embedded with smooth properly embedded boundary curves 
$\Lambda_i=(G\circ F)(C_i\setminus \{p_i\})$ $(i=1,\ldots,m)$ diffeomorphic to $\R$.
By using Theorem \ref{th:FR} and results of Stolzenberg \cite{Stolzenberg1966AM} it is then
possible to find a sequence of holomorphic automorphisms $\Phi_j\in \Aut(\C^2)$ converging on the interior 
$\mathring\Sigma = \Sigma\setminus \bigcup_{i=1}^m\Lambda_i$ of $\Sigma$ while the boundary 
curves $\Lambda_i$ diverge to infinity. If things are done right then the domain of convergence of 
the sequence $\Phi_j$ is a Fatou--Bieberbach domain $\Omega\subset\C^2$  
such that $\mathring \Sigma\subset \Omega\cong\C^2$
and $b\Sigma= \bigcup_{i=1}^m\Lambda_i \subset b\Omega$. This embeds the interior
$M\setminus bM$ of $M$ properly into $\Omega\cong\C^2$.

In Wold's joint paper with Forstneri\v c \cite{ForstnericWold2009} (2009), this construction was coupled 
with a newly developed technique of exposing boundary points of bordered Riemann surfaces. This led to  
the following result (see \cite[Theorem 1.1 and Corollary 1.2]{ForstnericWold2009}).
 
%
%   FW 2009
% 
\begin{theorem}\label{th:FW1}                  
Assume that $M$ is a compact bordered Riemann surface with $\Cscr^r$ boundary $(r>1)$. 
Every $\Cscr^1$ embedding $M\hra\C^2$ that is holomorphic in $\mathring M=M\setminus bM$
can be approximated uniformly on compacts in $\mathring M$ 
by proper holomorphic embeddings $\mathring M \hra\C^2$.          
\end{theorem}

This shows that, for bordered Riemann surfaces, the main part of Problem \ref{prob:RS}   
is to find a holomorphic embedding of the closed surface (including the boundary) in $\C^2$. 

In the same year 2009, Kutzschebauch, L{\o}w and Wold \cite{KutzschebauchLowWold2009}
provided examples of open Riemann surfaces $M$ which embed properly holomorphically
in $\C^2$ with interpolation, meaning that for every pair of discrete sequences $a_j\in M$ 
and $b_j\in\C^2$ without repetition there is a proper holomorphic embedding $F:M\hra \C^2$
with $F(a_j)=b_j$ for all $j=1,2,\ldots$.

It is natural to ask whether an analogue of Theorem \ref{th:FW1} also holds for Riemann surfaces
with infinitely many boundary curves. After some initial developments       
by Majcen \cite{Majcen2009,Majcen2013}, a fairly general result was obtained
by Forstneri\v c and Wold \cite[Theorem 5.1]{ForstnericWold2013} 
for domains in $\CP^1$ with at most countably many boundary components. 
By  He and Schramm \cite{HeSchramm1993}, 
such a domain is conformally equivalent to a  \emph{circled domain} $\Omega\subset\CP^1$, i.e.  
such that every connected component of $\CP^1 \setminus \Omega$ is a round disc or a point (puncture). 

%
%
%  FW 2013
%
%
\begin{theorem}\label{th:FW2}
Let $\Omega$ be a circled domain in $\CP^1$. If all but finitely many punctures in 
$\CP^1\setminus\Omega$ are limit points of discs in $\CP^1\setminus\Omega$, 
then $\Omega$ embeds properly holomorphically in $\C^2$.
\end{theorem}

The proof of this result in \cite{ForstnericWold2013} uses technical ingredients from the previous 
papers, but it relies on a considerably more delicate induction scheme. The problem is that the 
boundary components of $b\Omega$ may cluster on one another. The main new point is that 
at every step of the induction process one exposes and opens up a new Jordan curve in $b\Omega$ 
as described above, while at the same time pushes a carefully selected finite group of curves 
close to it towards infinity. The details are considerable.
Essentially the same proof gives the analogous result for circled domains in elliptic curves (tori). 

The problem of embedding Riemann surfaces with punctures properly into $\C^2$ is even more delicate,  
and no general techniques have been developed yet. Recently, 
Kutzsche\-bauch and Poloni \cite{KutzschebauchPoloni2020} (2020)
showed in particular that if $K$ is a countable closed subset of the Riemann sphere
$\CP^1$ with at most two accumulation points then the complement $\CP^1\setminus K$ 
admits a proper holomorphic embedding into $\C^2$. They proved the analogous result for 
complements of certain closed countable sets in complex tori and in hyperelliptic Riemann surfaces.

%
%
%   SECTION: LONG C^n's
%
%
\section{Complex manifolds exhausted by Euclidean spaces}\label{sec:long}

A complex manifold $X$ of dimension $n$ is called a {\em long $\C^n$}
if it is the union of an increasing sequence of domains 
$X_1\subset X_2\subset X_3\subset \cdots \subset \bigcup_{j=1}^\infty X_j=X$ such 
that each $X_j$ is biholomorphic to $\C^n$. By the Riemann mapping theorem, every long $\C$
is biholomorphic to $\C$. The first known example of a long $\C^n$ for any $n>1$ which is not
holomorphically convex, and hence not Stein, was found by Wold \cite{Wold2010} in 2010. 
In his example, every pair $X_k\subset X_{k+1}$ in the exhaustion corresponds to a non-Runge
Fatou--Bieberbach domain in $\C^n$ (see Theorem \ref{th:nonRungeFB}) such that, for some 
compact subset $K\subset X_1$, the $\Oscr(X_{k+1})$-hull of $K$ is not contained in $X_k$ 
for any $k\in\N$. It follows that the $\Oscr(X)$-hull of $K$ is not compact, hence $X$ is not Stein.
 
Subsequently, Forstneri\v c \cite{Forstneric2012PAMS} showed that for $n>1$ and any pair of 
disjoint countable sets $A,B\subset \C$ there is a holomorphic submersion $F:Z\to\C$ from 
an $(n+1)$-dimensional complex manifold $Z$ such that
\begin{enumerate}[\rm (a)]
\item every fibre $Z_z=F^{-1}(z)$ $(z\in \C)$ is a long $\C^n$, 
\item the fibre $Z_z$ is non-Stein for every $z\in A$, and
\item $Z_z$ is biholomorphic to $\C^n$ for every $z\in B$.
\end{enumerate}
By choosing $A$ and $B$ to be everywhere dense in $\C$, one gets a submersion onto $\C$
such that the type of the fibre jumps near every point of $\C$ from $\C^n$ to a non-Stein long $\C^n$. 
 
The questions whether there exist long $\C^2$'s without nonconstant holomorphic functions, 
or non-biholomorphic non-Stein long $\C^2$'s, were answered affirmatively by Boc Thaler
and Forstneri\v c \cite{BocThalerForstneric2016} in 2016. One of their results is the following.

%
%   LONG C^n WITHOUT NONCONSTANT FUNCTIONS
%
\begin{theorem}\label{th:noholomorphic}
For every $n>1$ there exists a long $\C^n$ without any nonconstant holomorphic 
or plurisubharmonic functions.
\end{theorem}

To prove this result, they used Wold's construction of 
a non-Runge Fatou--Bieberbach domain \cite{Wold2010}, 
but applied it inductively in a considerably more intricate manner.

Theorem \ref{th:noholomorphic} gives an essentially optimal counterexample to the classical 
{\em union problem for Stein manifolds}, asking whether an increasing union of Stein manifolds is always Stein. 
For domains in $\C^n$ this question was raised by Behnke and Thullen \cite{BehnkeThullen1934} in 1934, 
and an affirmative answer was given by Behnke and Stein \cite{BehnkeStein1939} in 1939.
Some progress on the general question was made by Stein \cite{Stein1956} and 
Docquier and Grauert \cite{DocquierGrauert1960}. The first counterexample to the union problem in 
any dimension $n\ge 3$ was found by Forn\ae ss \cite{Fornaess1976MA} in 1976.
He constructed an increasing union of balls that is not holomorphically convex, hence not Stein. 
His proof is based on an example of a biholomorphic map 
$\Phi: \Omega\stackrel{\cong}{\to} \Phi(\Omega)\subset \C^3$ 
from a bounded neighbourhood $\Omega\subset \C^3$ 
of any compact set $K\subset \C^3$ with nonempty interior 
such that the polynomial hull of $\Phi(K)$ is not contained in $\Phi(\Omega)$. 
(An example of this phenomenon was discovered by Wermer \cite{Wermer1959} already in 1959.)
In 1977, Forn{\ae}ss and Stout constructed a three-dimensional increasing union of polydiscs 
without any nonconstant holomorphic function \cite{FornaessStout1977}. 
Increasing unions of hyperbolic Stein manifolds were studied further by 
Forn{\ae}ss and Sibony \cite{FornaessSibony1981} and Forn{\ae}ss \cite{Fornaess2004}. 
For the connection with Bedford's conjecture we refer to the survey by Abbondandolo
et al. \cite{Abbondandolo2014}. 
In dimension $n=2$ the first counterexample to the union problem 
was the aforementioned example of Wold \cite{Wold2010} of a non-Stein long $\C^2$. 

To answer the second question concerning the existence of non-biholomorphic non-Stein 
long $\C^n$'s, Boc Thaler and Forstneri\v c introduced in \cite{BocThalerForstneric2016} 
new biholomorphic invariants of a complex manifold $X$, the {\em stable core} 
and the {\em strongly stable core}, which allow one to distinguish some long $\C^n$'s 
from one another. 

The stable core of $X$ is the set of all points $x\in X$ which 
admit a compact neighbourhood $K$ such that,
for some (and hence for any) increasing sequence of compact sets 
$K_1\subset K_2\subset \cdots \subset \bigcup_{j=1}^\infty K_j = X$ with $K\subset K_1$
and $K_j\subset \mathring K_{j+1}$ for all $j$, the increasing sequence of hulls
$\wh K_{\Oscr(K_j)}$ stabilizes at some $j=j_0\in\N$. (Such $K$ is said to have the
{\em stable hull property}). Hence, the stable core is an open subset of $X$. 
A compact set $K$ in $X$ is called the strongly stable core of $X$
if $K$ has the stable hull property but any compact set $L\subset X$ with
$\mathring L\setminus K \ne \varnothing$ fails to have the stable hull property.
In a Stein manifold or a compact manifold, the stable core is the entire manifold, 
and a Stein manifold does not have a strongly stable core. However, these invariants are often nontrivial 
in complex manifolds obtained as increasing unions of domains which fail to form Runge pairs. 
Boc Thaler and Forstneri\v c proved the following result 
\cite[Theorem 1.2]{BocThalerForstneric2016}. 

%
%
%  MANIFOLDS X(B)
%
%
\begin{theorem}\label{th:XB}
To every compact, strongly pseudoconvex and polynomially convex 
domain $B\subset \C^n$, $n>1$, we can associate a complex manifold $X(B)$, 
which is a long $\C^n$ containing a biholomorphic copy of $B$ as its strongly stable core, 
such that every biholomorphic map $\Phi : X(B)\to X(C)$ between two such
manifolds takes $B$ onto $C$. In particular, for every 
$\Phi\in\Aut(X(B))$ the restriction $\Phi|_B$ is a holomorphic automorphism of $B$. 
\end{theorem}

It follows that if $X(B)$ is biholomorphic to $X(C)$ then $B$ is biholomorphic to $C$. The 
construction likely gives many non-equivalent long $\C^n$'s associated to the same domain $B$. 

The manifold $X(B)$ is an increasing union 
$X_1\subset X_2\subset\cdots\subset \bigcup_{i=1}^\infty X_i=X$ of domains
$X_i\cong \C^n$ such that $B\subset X_1$, $\wh B_{\Oscr(X_i)}=B$ for every $i\in\N$,
but for any compact neighbourhood $K$ of a point $p\in X\setminus B$ 
the hull $\wh K_{\Oscr(X_{i+1})}$ is not contained in $X_i$ for large $i\in \N$.
Every inclusion $X_i\hra X_{i+1}$ is given by a Fatou--Bieberbach 
map $\phi_i:\C^n\hra \C^n$ such that $\phi_i(\C^n)$ is not Runge (like in Wold's
example \cite{Wold2008}), and 
the $\Oscr(X_{i+1})$-hull of the image of $K$ in $X_i$ intersects $X_{i+1}\setminus X_i$.
Thus, the increasing sequence of hulls $\wh K_{\Oscr(X_{i})}$ does not stabilize.

It was shown by Poincar\'e in 1907 \cite{Poincare1907} 
that most pairs of smooth strongly pseudoconvex hypersurfaces in $\C^n$ 
are not biholomorphic to each other.
A complete set of countably many local holomorphic invariants of such hypersurfaces 
is provided by the Chern--Moser normal form \cite{ChernMoser1974}. 
Hence, Theorem  \ref{th:XB} implies the following.

%
%
%  INFINITELY MANY LONG C^n's.
%
%
\begin{corollary}\label{cor:infinitely}
For every $n>1$ there exist uncountably many  non-equivalent long $\C^n$'s
such that none of them has any holomorphic automorphisms different from the identity map.
\end{corollary}

The following challenging problems remain open.

\begin{problem}
\begin{enumerate}[\rm (a)] 
\item Which open subsets $U\subset \C^n$ are the stable core of a long $\C^n$?
\item Is there a long $\C^n$ which is Stein but not biholomorphic to $\C^n$?
\item Is there a non-Stein long $\C^2$ with a nonconstant holomorphic function?
\item Is there a long $\C^2$ without any nonconstant meromorphic functions? 
\end{enumerate}
\end{problem}

Question (d) is motivated by the observation that some meromorphic functions survive
in the construction leading to Theorems \ref{th:noholomorphic} and \ref{th:XB}.

%
%
%   SECTION: STEIN MANIFOLDS WITH THE DENSITY PROPERTY AND OKA THEORY
%
%
\section{Stein manifolds with the density property and Oka manifolds}\label{sec:Oka}

In this section we discuss the role that Stein manifolds with the density property 
play in the theory of Oka manifolds, in holomorphic embedding problems, and in complex 
dynamics.

\smallskip
\noindent{\bf Stein manifolds with the density property are Oka.}
Recall (see \cite[Section 5.4]{Forstneric2017E}) that a complex manifold $X$ is said to be 
an {\em Oka manifold} if every holomorphic map $U\to X$ from a neighbourhood of any 
given compact convex set $K\subset \C^n$, $n\in\N$, can be approximated uniformly on 
$K$ by entire maps $\C^n\to X$. This {\em convex approximation property} (CAP) 
is one of several equivalent characterisations of Oka manifolds. 
Recall that holomorphic maps $S\to X$ from any 
Stein manifold $S$ to an Oka manifold $X$ satisfy all forms of the h-principle; 
see \cite[Theorem 5.4.4]{Forstneric2017E} for a precise statement.

Another more recent characterisation of Oka manifolds is due to Kusakabe \cite{Kusakabe2021IUMJ}. 
Consider the following condition on a complex manifold $X$:

(*) For any compact convex set $L\subset \C^m$ $(m\in\N)$, open set $U\subset \C^m$ containing $L$,  
and holomorphic map $f:U\to X$ there are an open set $V$ with $L\subset V\subset U$ and 
a holomorphic map $F:V\times \C^n\to X$ for some $n \ge \dim X$  such that $F(\cdotp,0)=f|_V$ and 
\[ %begin{equation}\label{eq:spray}
	\frac{\di}{\di t}\Big|_{t=0} F(z,t):\C^n \to T_{f(z)} X
	\ \ \text{is surjective for every}\ z\in V.
\] %end{equation}
Here, $t=(t_1,\ldots,t_n)\in\C^n$.  Such $F$ is called a {\em dominating holomorphic spray over $f|_V$}.

This is a restricted version of condition Ell$_1$ introduced by Gromov \cite[p.\ 72]{Gromov1986};
see also \cite{Gromov1989}. In \cite[Theorem 1.3]{Kusakabe2021IUMJ}, Kusakabe used the technique 
of gluing sprays from \cite[Sect.\ 5.9]{Forstneric2017E} to show that this condition implies CAP,
so a complex manifold satisfying (*) is an Oka manifold. Conversely, every Oka manifold satisfies 
condition $\mathrm{Ell}_1$ by \cite[Corollary 8.8.7]{Forstneric2017E}.

A holomorphic map $F:X\times \C^n\to X$ satisfying condition (*) with $V=X$ and $f=\Id_X$ 
is a {\em dominating spray on $X$}.
A complex manifold $X$ which admits a dominating spray is called {\em elliptic}. This terminology 
was introduced by Gromov \cite{Gromov1989}, who proved that every elliptic manifold is an Oka manifold; 
the details can be found in \cite{ForstnericPrezelj2000} and \cite[Chapter 5]{Forstneric2017E}. 
Conversely, every Stein Oka manifold is elliptic; see \cite[3.2.A]{Gromov1989} or
\cite[Proposition 5.6.15]{Forstneric2017E}.

Suppose now that $X$ is a Stein manifold with the density property. It is easily seen that 
complete holomorphic vector fields on $X$ span the tangent space $T_x X$ at every point,
so  $X$ is holomorphically  flexible in the sense of Arzhantsev et al.\  
\cite{ArzhantsevFlennerKalimanKutzschebauchZaidenberg2013DMJ}. 
It follows that for any compact subset $K$ of $X$ there exist finitely many complete 
holomorphic vector fields which generate $TX$ over $K$. Denoting their flows by $\phi^1,\ldots,\phi^n$, 
we obtain a spray $F:X\times \C^n \to X$ of the form 
\begin{equation}\label{eq:flowspray}
	F(x,t_1,\ldots,t_n)=\phi^1_{t_1}\circ \phi^2_{t_2}\circ\cdots\circ \phi^n_{t_n}(x),
	\quad x\in X,\ t_j\in \C, \ j=1,\ldots,n
\end{equation}
which is dominating at every point of $K$.
This means that $X$ is {\em weakly elliptic}, and hence an Oka manifold 
\cite[Corollary 5.5.12]{Forstneric2011E}.
%(Another argument is that flexibility evidently implies Condition $\mathrm{Ell}_1$ (see (*) above), and hence the manifold is Oka by Kusakabe's \cite[Theorem 1.3]{Kusakabe2021IUMJ}.) 
With some more work, one can show that finitely many complete holomorphic 
vector fields span the tangent bundle of $X$ at every point, and hence $X$ admits a 
globally dominating spray of the form \eqref{eq:flowspray}; see
\cite[Theorem 4]{KalimanKutzschebauch2008MZ} or 
\cite[Proposition 5.6.22 (b)]{Forstneric2017E}. Let us record this and a few related results. 

%
%   DENSITY IMPLIES OKA
%
\begin{theorem}\label{th:DPOka}
\begin{enumerate}[\rm (a)]
\item Every Stein manifold with the density property is an Oka manifold.
\item Every Stein manifold with the volume density property
% (with respect to a holomorphic volume form) 
is an Oka manifold.
\item Every complex manifold with the density property 
whose tangent bundle is pointwise spanned by globally defined 
holomorphic vector fields is an Oka manifold.
\end{enumerate}
\end{theorem}

The argument leading to statement (a) was given above (see also \cite[Theorem 5.5.18]{Forstneric2011E}).
Part (b) follows from \cite[Lemma 4.1]{KalimanKutzschebauch2011}, 
which says that on a Stein manifold with the volume density property there are finitely 
many complete divergence-free holomorphic vector fields which span the tangent space at each point,
so the manifold is elliptic. The main point is to use the isomorphism \eqref{eq:Theta} 
between holomorphic vector fields with vanishing divergence 
and closed holomorphic $(n-1)$-forms, where $n$ is the dimension of the manifold, and the fact
that Cartan's Theorem A provides many exact (hence closed) holomorphic $(n-1)$-forms. 
Part (c), which does not use Steinness, is obtained by noting 
that the two conditions  ensure that the tangent bundle $TX$ is spanned
over any compact subset of $X$ by finitely many $\C$-complete holomorphic 
vector fields, so $X$ is weakly elliptic and hence Oka.

%
%   Remark concerning Proposition 5.6.23 in \cite{Forstneric2017E}
%
\begin{remark}\label{rem:gap}
Theorem \ref{th:DPOka} (c) corresponds to \cite[Proposition 5.6.23]{Forstneric2017E}. 
In the latter source the hypothesis that the tangent bundle $TX$ is  
spanned by globally defined holomorphic vector fields is accidentally missing, 
but is tacitly used in the proof.
\end{remark}

The connection between Stein manifolds with the density property and Oka manifolds 
goes way beyond what has been said so far.
As shown in \cite{ForstnericWold2020MRL}, Theorem \ref{th:FW2020} together with the   
characterisation of Oka manifolds by Condition Ell$_1$ (see (*) above) easily implies the following 
result of Kusakabe \cite[Theorem 1.2 and Corollary 1.3]{Kusakabe2020complements} from 2020. 

%
%   KUSAKABE'S THEOREM
%
\begin{theorem} \label{th:KOka}
Let $X$ be a Stein manifold with the density property. For every compact $\Oscr(X)$-convex subset 
$K$ of $X$ the complement $X\setminus K$ is an Oka manifold. In particular, the complement 
of any compact polynomially convex set  in $\C^n$, $n>1$, is an Oka manifold. 
\end{theorem}

For $X=\C^n$ the proof goes as follows. 
Given a holomorphic map $f:V\to \C^n\setminus K$ as in (*) on a neighbourhood 
of a compact convex set $L\subset\C^m$, Theorem \ref{th:FW2020} gives 
a dominating spray $F:V\times \C^n \to \C^n\setminus K$
such that for every $z\in V$, $F(z,\cdotp):\C^n\to \C^n\setminus K$ is a 
Fatou--Bieberbach map and $F(z,0)=f(z)$. Thus,   
$\C^n\setminus K$ satisfies Condition $\mathrm{Ell}_1$,  
and hence is Oka by Kusakabe's theorem \cite{Kusakabe2021IUMJ}.
The general case follows from \cite[Theorem 3.1]{ForstnericWold2020MRL}. 

Kusakabe also proved that complements of certain unbounded 
closed polynomially convex sets $E\subset\C^n$ are Oka. (A closed set 
is polynomially convex if it is exhausted by an increasing sequence 
of compact polynomially convex sets.) The following is
\cite[Theorem 1.6]{Kusakabe2020complements}. 

%
%   KUSAKABE: complements of unbounded polynomially convex sets
%
\begin{theorem}\label{th:E}
Let $E$ be a closed polynomially convex set in $\C^n$, $n\ge 2$, such that
\begin{equation}\label{eq:E}
	E \subset \bigl\{(z,w)\in \C^{n-2}\times \C^2: |w|\le C(1+|z|)\bigr\}
\end{equation}
for some $C>0$. Then $\C^n\setminus E$ is an Oka manifold.
\end{theorem}

To prove this result, Kusakabe first showed that %, under the hypotheses of the theorem, 
the restricted coordinate projection $\pi:\C^n\setminus E\to \C^{n-2}$, $\pi(z,w)=z$, has the
Oka property for liftings; see \cite[Corollary 5.5.11]{Forstneric2017E} for this notion.
In particular, any holomorphic map $f:V\to \C^n\setminus E$ from a Stein manifold $V$ 
(again, it suffices to consider convex domains in Euclidean spaces)
is the core of a fibre-dominating spray $F:V\times \C^m \to \C^n\setminus E$
with $F(\cdotp, 0)=f$ and $\pi\circ F(z,t)=\pi\circ f(z)$ for all $z\in V$ and $t\in\C^m$.
By \eqref{eq:E} the same holds for linear projections 
$\pi':\C^n \to \C^{n-2}$ with kernels close to $\ker\pi=\{0\}^{n-2}\times \C^2$.
This gives sprays over $f$ with values in $\C^n\setminus E$ 
which are fibre-dominating in directions spanning the tangent space
to $\C^n$ at any point. By \cite[Corollary 4.1]{Kusakabe2021IUMJ} it follows
that $\C^n\setminus E$ is Oka. (The argument amounts to composing such sprays
to obtains a dominating spray over $f$ and then applying 
\cite[Theorem 1.3]{Kusakabe2021IUMJ}.)

Kusakabe also proved the following result (see \cite[Theorem 4.2]{Kusakabe2020complements}),
which is an interesting and powerful application of the fibred density property discussed
in Subsection \ref{ss:fibreddensity}. 

%
%  KUSAKABE: THEOREM 4.2 (2020)
%
\begin{theorem}\label{th:Kusakabe-fibreddensity}
Let $\pi : Y\to B$ be a holomorphic submersion between reduced
complex spaces. Assume that $E$ is a closed subset of $Y$ such that
every point $b\in B$ admits an open neighbourhood $U\subset B$ satisfying
the following conditions:
\begin{enumerate}[\rm (i)]
\item $Y_U:=\pi^{-1}(U)$ is a Stein space,
\item $E_U:=E\cap \pi^{-1}(U)$ is holomorphically convex in $Y_U$, and 
\item the projection $\pi: Y_U\setminus E_U\to U$ enjoys the fibred density property
(see Subsec.\ \ref{ss:fibreddensity}).
\end{enumerate}
Then the restriction $\pi:Y\setminus E\to B$ enjoys the Oka property for lifings.
\end{theorem}

We also mention the result of Forstneri\v c and L\'arusson \cite{ForstnericLarusson2014}
which says that the holomorphic automorphism group $\Aut(\C^n)$ for $n\ge 2$ enjoys most Oka 
properties for holomorphic maps $X\to \Aut(\C^n)$ from Stein manifolds.
Although $\Aut(\C^n)$ does not carry the structure of an (infinite dimensional) complex
manifold, it is natural to consider a map $f:X\to \Aut(\C^n)$ holomorphic if the associated 
evaluation map $F:X\times\C^n \to \C^n$, given by $F(x,z)=f(x)(z)$ $(x\in X,\ z\in\C^n)$, is holomorphic.
It is an open problem whether the same is true for the automorphism group of every 
Stein manifold with the density property in place of $\C^n$.
In fact, it seems that this is not known for any other example besides $\C^n$.

\smallskip
\noindent{\bf Embedding Stein manifolds in Stein manifolds with the density property.}
%
%   EMBEDDINGS INTO STEIN MANIFOLDS WITH DENSITY PROPERTY
%
We consider the problem of embedding Stein manifolds into model complex 
manifolds. The ideal class of models would be Oka manifolds. 
The following result (see \cite[Corollary 8.9.3]{Forstneric2017E}) 
is obtained by combining the main result of Oka theory 
\cite[Theorem 5.4.4]{Forstneric2017E} with the jet transversality theorem for holomorphic maps 
from Stein manifolds to Oka manifolds \cite[Theorem 8.9.1]{Forstneric2017E}. 

%
%   Injective immersions
%
\begin{theorem}
Let $X$ be a Stein manifold and $Y$ be an Oka manifold. If $\dim Y\ge 2\dim X$,
then every continuous map $X\to Y$ is homotopic to a holomorphic immersion with simple 
double points. If $\dim Y\ge 2\dim X+1$ then the immersion can be chosen injective.
\end{theorem}

On the other hand, there are Oka manifolds which do not admit any {\em proper} 
holomorphic images of even the simplest manifolds such as the disc; 
see \cite[Example 1.3]{DrinovecForstneric2010AJM} which involves certain 
punctured tori of dimension $>1$ (these are Oka). 
Stein manifolds with the density property are much better in this respect, as 
demonstrated by the following result.

%
%  EMBEDDING-IMMERSION THEOREM
%
\begin{theorem}\label{th:EmbDP}
Let $X$ and $Y$ be Stein manifolds, and assume that $Y$ has the density property
or the volume density property. Then, the following hold.
\begin{enumerate}[\rm (a)]
\item If $2\dim X +1\le \dim Y$ then any continuous map $f: X\to Y$ is homotopic to a proper holomorphic 
embedding $F: X\hra Y$. If in addition $f$ is holomorphic on a neighborhood of a compact $\Oscr(X)$-convex 
set $K\subset X$ and $X'$ is a closed complex subvariety of $X$
such that the restriction $f|_{X'} : X'\hra Y$ is a proper holomorphic embedding, 
then $F$ can be chosen to agree with $f$ on $X'$ and to approximate $f$ uniformly on $K$.
\item 
If $2\dim X=\dim Y$ then any continuous map $X\to Y$ is homotopic to a proper holomorphic 
immersion $X\to Y$ with simple double points, with additions as in part (a) concerning 
approximation and interpolation.
\end{enumerate}
\end{theorem}

Part (a) was proved by Andrist, Forstneri\v c, Ritter, and Wold \cite{AndristForstnericRitterWold2016}
in 2016. (The special case for Riemann surfaces was obtained beforehand by 
Andrist and Wold \cite{AndristWold2014}.) Part (b) is due to 
Forstneri\v c \cite{Forstneric2019JAM} (2019). 
The proofs strongly depend on the Anders\'en--Lempert theory. Using Theorem \ref{th:FR2}, 
one inductively constructs a sequence of continuous maps $X\to Y$ which are 
holomorphic embeddings (or immersions) on larger and larger domains in $X$ such that 
the sequence converges uniformly on compacts 
to a proper holomorphic embedding or immersion $X\to Y$. The assumption that $Y$ 
has the (volume) density property is crucial in this proof.

Theorem \ref{th:EmbDP} is classical when $Y$ is a Euclidean space $\C^N$ \cite{Remmert1956,Narasimhan1960AJM,Bishop1961AJM}.  
In this case, the optimal embedding dimension is $N=\left[\frac{3\dim X}{2}\right] +1$
if $\dim X>1$ according to Eliashberg and Gromov \cite{EliashbergGromov1992AM} and 
Sch\"urmann \cite{Schurmann1997} (see also Section \ref{sec:RS}).
It is not known whether the embedding or immersion dimension in Theorem \ref{th:EmbDP} 
can be lowered for more general Stein manifolds with the density property as targets.
Another problem is the following.

\begin{problem} 
Does Theorem \ref{th:EmbDP} hold for every Oka Stein manifold $Y$?
\end{problem}

The following is a corollary to Theorem \ref{th:EmbDP} (b) and the fact that 
the space $(\C^*)^n$ with coordinates $z=(z_1,\ldots,z_n)$ %(where $z_j\in \C^*$ for $j=1,\ldots,n$)
enjoys the volume density property with respect to the volume form 
$
	\omega= \frac{dz_1\wedge\cdots\wedge dz_n}{z_1\cdots z_n}.
$
(See \cite{Varolin2001} or \cite[Theorem 4.10.9 (c)]{Forstneric2017E}.) 

\begin{corollary}\label{cor:harmonic}
Every Stein manifold $X$ of complex dimension $n\ge 1$ admits a proper holomorphic immersion 
to $(\C^*)^{2n}$ and a proper pluriharmonic map to $\R^{2n}$.
\end{corollary}

This provides a counterexample in any dimension to the conjecture of  Schoen and Yau \cite{SchoenYau1997} 
that the unit disc $\D$ does not admit any proper harmonic maps to $\R^2$. We refer to 
\cite[Subsect.\ 3.3]{Forstneric2018Korea} and \cite[Sect.\ 3.10]{AlarconForstnericLopez2021} 
for the discussion of this topic, also in the context of minimal surfaces,  
and for references to the earlier counterexamples for the disc $X=\D$.

\smallskip
\noindent{\bf Open Runge embeddings.} 
%
%   RUNGE TUBES IN STEIN MANIFOLDS WITH DP
%
%
%We mention another circle of results in the same general framework.
%They concerns open holomorphic embeddings of certain Stein manifolds onto 
%Runge domains in Stein manifolds with the density property.
%
It was a long-standing problem whether  $\C^*\times\C$ embeds 
as a Runge domain in $\C^2$; such hypothetical domains were called {\em Runge cylinders} 
in $\C^2$. This question arose in connection with the classification of Fatou components for 
H\'{e}non maps by Bedford and Smillie \cite{BedfordSmillie1991JAMS} in 1991. 
In 2021, Bracci, Raissy and Stens\o nes \cite{BracciRaissyStensones2021} 
obtained a Runge embedding $\mathbb C^*\times\mathbb C\hra \mathbb C^2$
as the basin of a non-polynomial holomorphic automorphism of $\C^2$
at a parabolic fixed point. In 2020, Forstneri\v c and Wold \cite{ForstnericWold2020PAMS} showed that 
Runge tubes are abundant in Stein manifolds with the density property. 
Although their proof is completely different from 
the one in \cite{BracciRaissyStensones2021}, both use Anders\'{e}n-Lempert theory. 
Their first result \cite[Theorem 1.1]{ForstnericWold2020PAMS} is the following.

%
%   FW2020: MAIN THEOREM 1
%
\begin{theorem}\label{th:FW2020Main1}
Let $X$ and $Y$ be Stein manifolds with $\dim X<\dim Y$, and assume that $Y$ has
the density property.  Suppose that $\theta:X\hra Y$ is a holomorphic embedding with 
$\Ocal(Y)$-convex image (this holds in particular if $\theta$ is proper), and let $E\to X$ 
denote the holomorphic normal bundle associated to $\theta$. Then, $\theta$ can be approximated 
uniformly on compacts in $X$ by  holomorphic embeddings of $E$ into $Y$ whose 
images are Runge domains in $Y$.  
\end{theorem}

To get a Runge embedding of $\mathbb C^*\times\mathbb C$ into $\mathbb C^2$,
it suffices to embed $X=\mathbb C^*$  onto the curve $\{zw=1\}\subset \mathbb C^2$ and note 
that any holomorphic vector bundle over $\mathbb C^*$ (indeed, over any open Riemann surface) 
is trivial by Oka's theorem \cite{Oka1939}. (See also \cite[Sect.\ 5.2]{Forstneric2017E}.)
This argument gives the following corollary to Theorem \ref{th:FW2020Main1}.
%, noting also that every open Riemann surface admits a proper holomorphic embedding into $\C^3$.

%
%  COROLLARY: RUNGE TUBES OVER OPEN RIEMANN SURFACES
%
\begin{corollary}[Runge tubes over open Riemann surfaces] \label{cor:tubesC2} 
If $X$ is an open Riemann surface which admits a proper holomorphic embedding into $\C^2$,
then $X\times \C$ is biholomorphic to a Runge domain in $\C^2$. For every open Riemann surface $X$
and $k\ge 2$, $X\times\C^k$ admits a Runge embedding into $\C^{k+1}$, and 
into any Stein manifold $Y^{k+1}$ with the density property.
\end{corollary}

The Runge embedding $E\hra Y$ of the normal bundle in Theorem \ref{th:FW2020Main1} 
need not agree with the given embedding $\theta:X\hra Y$ on the zero section $X$ of $E$; 
this is impossible in general due to Theorem \ref{th:hyperboliccomplement}. 
However, we can ensure this condition for algebraic embeddings of codimension at least 2 into 
$\C^n$; see \cite[Theorem 1.4 and Corollary 1.5]{ForstnericWold2020PAMS}.

%
%   RUNGE TUBES: MAIN THEOREM 2
%
\begin{theorem}\label{th:FWmain2}
Let $X$ be a Stein manifold and $\theta: X\hra \C^n$ be a proper holomorphic embedding onto an
algebraic submanifold of $\C^n$. If $n\ge \dim X+2$ then $\theta$ extends to a 
holomorphic Runge embedding $E\hra \C^n$ of the total space of the normal bundle of $\theta$. 
In particular, every proper algebraic embedding $X\hra\C^{n}$ $(n\ge 3)$ of an affine algebraic curve 
extends to a holomorphic Runge embedding $X\times \C^{n-1}\hra \C^{n}$. 
\end{theorem}

The proof of Theorem \ref{th:FWmain2} uses the result of 
Kaliman and Kutzschebauch \cite[Theorems 4 and 6]{KalimanKutzschebauch2008IM} 
that the Lie algebra $\Liealg(\C^n, Y)$ of algebraic vector fields on $\C^n$
vanishing on an algebraic submanifold $Y\subset\C^n$ of dimension $\dim Y\le n-2$ 
enjoys the strong algebraic density property (see Example (1) in Subsection \ref{ss:relativedensity}).  

\smallskip
\noindent{\bf Holomorphic dynamics on Stein manifolds with the density property.}
%
%   AROSIO AND LARUSSON: DYNAMICS ON STEIN MANIFOLDS WITH DP
%
Holomorphic dynamics is a lively field of complex analysis. So far, dynamical phenomena have mainly
been studied on complex Euclidean spaces, which have an abundance of holomorphic endomorphisms
and automorphisms. It is natural to ask to what extent can these results be extended to Stein 
manifolds with the density property, and how does the possibly nontrivial topology
of the manifold impact the global behaviour of maps under consideration.

The first steps in this emerging field were made recently by Arosio and L\'arusson. 
In \cite{ArosioLarusson2018} they proved that automorphisms with chaotic behaviour 
are generic among volume preserving automorphisms of a Stein manifold $X$ 
having the density property for an exact volume form. For $X=\C^n$, $n\ge 2$, 
with the standard volume form this was proved by Forn\ae ss and Sibony \cite{FornaessSibony1997},
and the authors follow their approach. 
They also showed that a generic volume preserving automorphism 
has a hyperbolic fixed point whose stable manifold is dense in $X$, generalizing 
a result of Peters, Vivas, and Wold on $\C^n$ \cite{PetersVivasWold2008}. 
In their second paper \cite{ArosioLarusson2019}, they proved closing lemmas for automorphisms of 
a Stein manifold with the density property and for endomorphisms of an Oka Stein manifold. 
In the former case they needed to impose a new tameness condition. It follows that hyperbolic 
periodic points are dense in the tame non-wandering set of a generic automorphism of a Stein manifold 
with the density property and in the non-wandering set of a generic endomorphism of an Oka 
Stein manifold.

%
%
%   SECTION: COMPLETE COMPLEX SUBMANIFOLDS
%
%
\section{Complete complex submanifolds}\label{sec:complete}

In 1977, Paul Yang asked \cite{Yang1977AMS,Yang1977JDG} whether 
there exist bounded complete immersed or embedded complex submanifolds in 
complex Euclidean spaces. Here, an immersion $X\to\C^n$ is said to be complete 
if the image of any divergent path $\gamma:[0,1)\to X$ (i.e., one that leaves
any compact subset of $X$ as $t\to 1$) has infinite Euclidean length. 
Equivalently, the Riemannian metric on $X$ induced by the immersion is a complete metric.
Yang's problem is an analogue of the {\em Calabi-Yau problem} for 
minimal surfaces in $\R^n$. We refer to \cite[Chapter 7]{AlarconForstnericLopez2021}
for background and a discussion of recent results on this subject.
 
The first such examples were constructed by Jones in 1979 \cite{Jones1979}
(a immersed holomorphic disc in the ball $\B^2$ of $\C^2$, and an embedded one 
in the ball of $\C^3$). Much later, Alarc\'on and Forstneri\v c
proved in \cite{AlarconForstneric2013MA} (2013) that every finite bordered Riemann surface admits
a complete proper holomorphic immersion in $\B^2$ and embedding in $\B^3$.
For further and more recent developments on this subject we refer to the papers 
\cite{AlarconJDG,AlarconForstneric2020MZ,AlarconGlobevnik2017,AlarconGlobevnikLopez2019Crelle,AlarconLopez2016JEMS,Drinovec2015JMAA,Globevnik2015AM}.
% and the 2019 survey \cite{AlarconForstneric2019JAMS} by Alarc\'on and Forstneri\v c. 
In particular, the following is a compilation of results by 
Globevnik \cite{Globevnik2015AM,Globevnik2016MA}, Alarc\'on, Globevnik and L\'opez
\cite{AlarconGlobevnikLopez2019Crelle}, Alarc\'on \cite[Corollary 1.2]{AlarconJDG}, and Alarc\'on and 
Forstneri\v c \cite{AlarconForstneric2020MZ}.  

%
%  COMPLETE COMPLEX CURVES IN THE BALL
%
\begin{theorem}\label{th:complete} 
For every pair of integers $1\le q<n$ there exists a holomorphic submersion 
$f:\B^n\to \C^q$ whose fibres are complete complex submanifolds of $\B^n$. 
In particular, the ball $\B^n$ can be foliated by complete properly embedded holomorphic discs.
\end{theorem}

The techniques used in the last three mentioned papers rely on Anders\'en--Lempert theory,
and they provide some additional information on topology of the leaves. 
The main idea is to place in $\B^n$ a polynomially convex labyrinth consisting of countably 
many pairwise disjoint compact sets (the authors used closed balls in affine hyperplanes) such that 
any divergent curve in $\B^n$ avoiding all but finitely many pieces 
of the labyrinth has infinite length. One then uses holomorphic automorphisms of $\C^n$ 
to successively twist a given holomorphic foliation on $\C^n$ such that the resulting sequence of
foliations converges in $\B^n$ to a foliation each of whose leaves 
avoids all but finitely many pieces of the labyrinth, so it is complete. 

On the other hand, Globevnik's technique in \cite{Globevnik2015AM,Globevnik2016MA} 
relies on a different idea, based on the construction of holomorphic functions on the ball 
which grow sufficiently fast on pieces of a suitable labyrinth. This technique does not allow 
any control of the topology of the leaves.

The methods in the cited papers, together with the construction of suitable labyrinths
by Charpentier and Kosi\'nski  \cite{CharpentierKosinski2020} (2020), 
show that the same result holds in all pseudoconvex Runge domains in $\C^n$, $n>1$.
Globevnik extended his original construction of complete complex submanifolds of $\B^n$ 
\cite{Globevnik2015AM} to this more general setting in \cite{Globevnik2016MA} (2016), 
but in his case (for constructing fast growing holomorphic functions)
less precise labyrinths suffice. 

In a related direction, Alarc\'on and Forstneri\v c constructed in \cite{AlarconForstneric2018PAMS} 
a complete injective holomorphic immersion $\C\to\C^2$ whose image is dense in $\C^2$. 
The analogous result was obtained for any closed complex submanifold $X\subset \C^n$
$(n>1)$ in place of $\C$. 

%
%
%	SECTION: AN APPLICATION TO 3-DIMENSIONAL TOPOLOGY
%
%
\section{An application in 3-dimensional topology}\label{sec:3D}

It is well-known that diffeomorphism groups of smooth or real analytic manifolds of positive dimension 
are huge, in particular infinite dimensional. Complete vector fields on such manifolds 
can be constructed using cutoff functions and approximation in the strong Whitney topology
(to obtain real analytic ones). It seems therefore not very interesting to study real shears on $\R^n$ 
or, say, multiplicative shears on a torus. 

However, there are situations where real shears gain importance for some specific reason. 
In this section we describe a recent example of this type.
It concerns the affirmative answer, given by Rafael Zentner \cite{Zentner2018} in 2018,
to the long standing problem in 3-dimensional topology asking whether the fundamental group of any 
homology 3-sphere different from the 3-sphere $S^3$ admits an irreducible representation into 
$\mathrm{SL}_2 (\C)$, i.e. a 2-dimensional irreducible representation. Here, a homology 3-sphere 
is a compact 3-manifold $X$ whose homology groups are those of $S^3$. 
Its fundamental group is nontrivial unless $X$ is $S^3$ itself.

%The affirmative answer to this problem, given by Rafael Zentner \cite{Zentner2018} in 2018,
%is an example where shears play a role in the real setting.  

Let us explain the main points in Zentner's analysis. 
There are three types of homology $3$-spheres, and the one type where the answer was not known 
are those homology $3$-spheres which admit a degree $1$-map to a splicing of two nontrivial knots in $S^3$. 
Since a degree $1$-map is $\pi_1$-surjective, 
it remained to find the answer for those homology $3$-spheres which arise by the 
following construction, called splicing. Take a pair of knots $K_1$, $K_2$ in $S^3$. Remove the 
tubular neighborhood $N(K_i)$ of the knot from $S^3$, i.e. set $X_i := S^3 \setminus  N(K_i)$, 
and glue these two manifolds along the boundaries of $N(K_i)$ (each isomorphic to a 2-dimensional torus) 
so that the longitude of one torus is glued onto the meridian of the other torus. The resulting manifold  
is a homology 3-sphere, called the {\em splicing} $Y_{K_1, K_2}$ of the knots $K_1$ and $K_2$. 

In order to produce an irreducible representation of the fundamental group 
$\pi_1 (Y_{K_1, K_2})$ into $\mathrm{SL}_2 (\C)$, it was important
to prove that for any two nontrivial knots $K_1$ and $K_2$ the images of their representation varieties in 
the representation variety $R(T^2)$ of the boundary torus $T^2$  (the place where we glued) intersect 
and thus yields a representation of $\pi_1 (Y_{K_1, K_2})$. 
This is the place where the real version of the Anders\'en-Lempert theory comes into play.

\begin{figure}[h!]
\def\svgwidth{0.6\columnwidth}
\begingroup%
  \makeatletter%
  \providecommand\color[2][]{%
    \errmessage{(Inkscape) Color is used for the text in Inkscape, but the package 'color.sty' is not loaded}%
    \renewcommand\color[2][]{}%
  }%
  \providecommand\transparent[1]{%
    \errmessage{(Inkscape) Transparency is used (non-zero) for the text in Inkscape, but the package 'transparent.sty' is not loaded}%
    \renewcommand\transparent[1]{}%
  }%
  \providecommand\rotatebox[2]{#2}%
  \ifx\svgwidth\undefined%
    \setlength{\unitlength}{368.50393066bp}%
    \ifx\svgscale\undefined%
      \relax%
    \else%
      \setlength{\unitlength}{\unitlength * \real{\svgscale}}%
    \fi%
  \else%
    \setlength{\unitlength}{\svgwidth}%
  \fi%
  \global\let\svgwidth\undefined%
  \global\let\svgscale\undefined%
  \makeatother%
  \begin{picture}(1,0.61538463)%
    \put(0,0){\includegraphics[width=\unitlength]{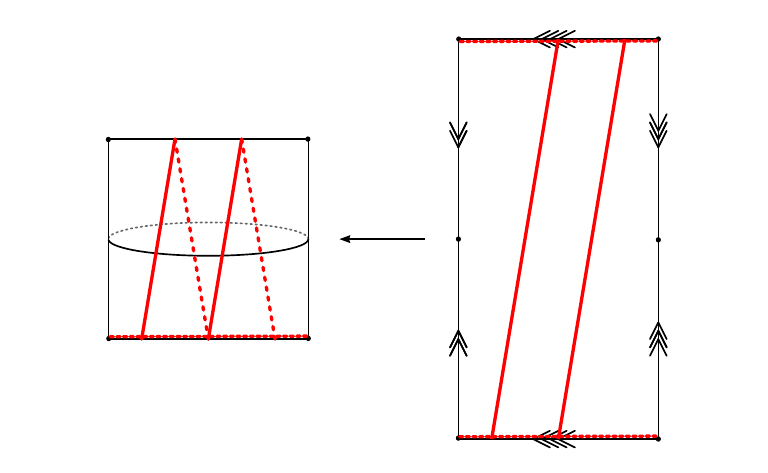}}%
    \put(0.3300436,0.38135272){\color[rgb]{1,0,0}\makebox(0,0)[lb]{\smash{$i^*(R(K))$}}}%
    \put(0.75062192,0.11094933){\color[rgb]{1,0,0}\makebox(0,0)[lb]{\smash{$i^*(R(K))$}}}%
  \end{picture}%
\endgroup%
\caption{The gluing pattern for obtaining the pillowcase from a rectangle, 
and the image of the representation variety $R(K)$ of the trefoil knot in the pillowcase}
\label{pillowcase trefoil}
\end{figure}

The representation variety, i.e. the space of $SU(2)$-representations of the fundamental group 
of a two-dimensional torus $T^2$ modulo conjugation,
\[
	R(T^2) = \Hom(\Z^2, SU(2))/SU(2)
\]
is homeomorphic to the {\em pillowcase}, a 2-dimensional sphere. 
In fact, if we denote generators of $\pi_1(T^2) \cong \Z^2$ by $m$ and $l$, 
then for a representation $\rho$ we may suppose that  
\[
	\rho(m) = \begin{bmatrix} e^{i \alpha} & 0 \\ 0 & e^{-i \alpha} \end{bmatrix} \hspace{0.5cm} 
	\text{and} \hspace{0.5cm}
	\rho(l) = \begin{bmatrix} e^{i \beta} & 0 \\ 0 & e^{-i \beta} \end{bmatrix},
\]
and hence we can associate to $\rho$ a pair $(\alpha,\beta) \in [0,2\pi] \times [0,2 \pi]$, which we also can 
think of as being a point on the two-dimensional torus $T ^2= \R^2 / 2 \pi \Z^2$. However, it is easily seen 
that a representation to which we associate $(2\pi - \alpha, 2 \pi - \beta)$ is conjugate to $\rho$. This is 
the only ambiguity, as the trace of an element in $SU(2)$ determines its conjugacy class. 
Therefore $R(T^2)$ is isomorphic to the quotient of the torus $T^2$ by the hyperelliptic involution 
$\tau: (\alpha,\beta) \mapsto (-\alpha,-\beta)$. This has four fixed points, and its quotient
\begin{equation*}
		R(T^2) = T/\tau
\end{equation*} 
is homeomorphic to the 2-sphere. It can also be seen as the quotient of the fundamental domain 
$[0,\pi] \times [0,2 \pi]$ for $\tau$ by identifications on the boundary as indicated in 
Figure \ref{pillowcase trefoil}. 

All abelian representations map to the red line `at the bottom'.

As explained above, one has to understand the image (under the map induced by restriction) of the 
representation variety $R(K)$ of the knot complement $S^3 \setminus K$ in the representation variety of 
the boundary torus $R(T^2)$ (the pillowcase). As an example, in Figure \ref{pillowcase trefoil} the image 
of $R(K)$ for $K$ being the trefoil knot is drawn in red.

Using results of Kronheimer and Mrowka, for any nontrivial knot $K$ the existence of points in the image 
of $R(K)$ lying on the straight line `on the top' of the pillowcase can be deduced, as well as on any 
holonomy perturbation of this top line.

The simplest class of holonomy perturbations are shearing maps on the torus,
\[
	\Phi_f : \R^2 / 2\pi \Z^2 \to \R^2 / 2\pi \Z^2, \quad (x,y) \mapsto (x, y + f(x))
\]
for some odd $2\pi$-periodic function $f$. Of course one can also change the roles of $x$ and $y$.
The following version of Theorem \ref{th:FR} for a real torus is due to 
Zentner \cite[Theorem 3.3]{Zentner2018}.

%
% SHEARS ON A TORUS
%
\begin{theorem}\label{th:torusshear} 
Let $\psi :  T \to T $ be an area-preserving map of the $n$-dimensional torus $T$ for $n\ge 2$,
which is volume preservingly isotopic to the identity, and let $\epsilon > 0$ be given. Then
there is a finite composition of shearing maps $\phi : T \to T$ which is $\epsilon$-close to $\psi$.
Moreover, the whole isotopy can be realized $\epsilon$-close to an isotopy through finitely
many shearing isotopies.
\end{theorem}

This result, or rather an equivariant version of it with respect to the involution, was the crucial ingredient  
to prove the existence of points in the image of $R(K)$ in the pillowcase $R(T^2)$ on any 
path connecting the left and right upper corners.
(Namely, using Theorem \ref{th:torusshear} and the Moser trick, any path can be approximated arbitrarily 
well by the straight line on the top of the pillowcase moved by holonomy perturbations; here, 
 shears. A closedness argument gives the desired conclusion.)
This latter fact can in turn be used to show that the images of the  representation varieties of the two 
nontrivial knots in the pillowcase always intersect. Indeed, the image of the representation variety 
of $K_1$ has to wrap around the pillowcase in one direction (its image undern the lift to the torus contains 
a path in the homology class of the longitude), while the image of the representation variety 
of $K_2$ has to wrap around the pillowcase in the other direction (its image under the lift to the torus 
contains a path in the homology class of the meridian). Thus, these images meet, thereby 
yielding a nontrivial representation of $Y_{K_1, K_2}$.

%
%
%   SECTION: THE RECOGNITION PROBLEM
%
%
\section{The recognition problem for complex Euclidean spaces}\label{sec:recognition}

The complex affine spaces $\C^n$  are the most natural and basic objects in algebraic and complex analytic 
geometries, comparable to the role played by the real Euclidean spaces for topological or differentiable 
manifolds. One of the most famous achievements of geometric topology for open real manifolds is the 
complete solution of the Open Poincar\'e Conjecture, namely, the characterisation of Euclidean spaces 
among open topological manifolds as the unique ones which are simultaneously contractible and 
simply connected at infinity. This conjecture was established by Stallings 
for PL-manifolds of dimension $n\geq 5$, Freedman in the case $n=4$, 
and finally Perelman in dimension $n=3$.

Nothing even remotely close to this exists in the complex algebraic or holomorphic case,
with the exception of the work of Ramanujam \cite{Ramanujam} on affine algebraic surfaces
(see below). Indeed, even extremely simple affine algebraic manifolds, such as the Koras--Russell 
cubic threefold KR in $\C^4$ defined by the polynomial equation $x+x^2 y +s^2+t^3= 0$ (see \eqref{eq:KR}), 
are not fully understood in this context. 

In his landmark paper \cite{Ramanujam} (1971), Ramanujam 
proved that a smooth contractible affine algebraic surface is 
isomorphic to $\C^2$ if it is simply connected at infinity. At the same time, he constructed
many examples of smooth contractible affine algebraic surface with non-trivial fundamental group 
at infinity, thereby opening the way for the construction of higher dimensional algebraic or analytic 
varieties which were later on verified, thanks to the cancellation theorems due to Iitaka--Fujita in the 
algebraic case and Zaidenberg in the holomorphic one, to be diffeomorphic to Euclidean spaces, 
while neither algebraically nor holomorphically isomorphic to $\C^n$. In the holomorphic case, 
even dimension $2$ remains  a mystery. These types of varieties, nowadays called 
{\em exotic affine spaces}, are challenging objects. 

These examples show that additional properties must be imposed on a variety which is 
diffeomorphic to $\R^{2n}$ in order to be biholomorphic or algebraically isomorphic 
to $\C^n$. A natural attempt is to use symmetries.
Affine spaces are homogeneous under the action of algebraic or holomorphic one-parameter flows. 
Recall that vector fields generating algebraic flows are called 
{\em locally nilpotent derivations}, LND's for short. 
A seminal breakthrough of Makar--Limanov \cite{MakarLimanov1996} in 1996 was to realize that 
the Koras--Russell cubic KR has some rigidity with respect to such algebraic flows. 
He introduced an invariant, now called the {\em Makar-Limanov invariant}, which measures the 
richness of algebraic flows and the homogeneity of a variety under the action of the group 
which they generate. Makar--Limanov and Kaliman developed in \cite{KalimanMakarLimanov1997} 
sophisticated algebraic techniques to compute this invariant for special classes of
affine varieties containing the Koras--Russell cubic. 
Refinements of these techniques by Dubouloz, Moser, Jauslin and Poloni 
enabled the construction of pairs of exotic affine $3$-folds failing the 
Zariski Cancellation Problem (see Problem \ref{prob:ZCP}), 
the construction of holomorphically trivial deformations 
of pairwise non-isomorphic algebraic exotic affine 3-folds \cite{DuboulozMoserPoloni2010},
and a full description of the group of algebraic automorphisms of the Koras-Russell cubic KR
\cite{MoserJauslin2009,DuboulozMoserPoloni2014}. The latter group is infinite dimensional like the group
$\mathrm{Aut}_{\mathrm{alg}}(\C^3)$, % of algebraic automorphisms of $\C^3$ but, surprisingly, 
but it acts on KR with precisely four orbits, in contrast to the transitivity of the action of 
$\mathrm{Aut}_{\mathrm{alg}}(\C^3)$ on $\C^3$.

In general, a characterisation of affine spaces is considered very useful if it can be used 
to solve the following Zariski Cancellation Problem. 

%
%  ZARISKI CANCELLATION PROBLEM
%
\begin{problem}\label{prob:ZCP}
Assuming that $X$ is an affine algebraic variety such that $X\times \C^k$ is isomorphic to $\C^{n+k}$
for some $k\in\N$, is $X$ is isomorphic to the affine space $\C^n$? The same question for $X$ a 
Stein manifold and the word isomorphic replaced by biholomorphic.
\end{problem}

For $n=1$ the affirmative answer is not difficult, while for $n=2$ it is a deep result due to 
Fujita \cite{Fujita1979} and Miyanishi and Sugie \cite{MiyanishiSugie1980}.
The case $n\ge 3$ is completely open. In the holomorphic case dimension 
$n=1$ is easy, and the problem is open in higher dimensions.

Since the density property of a Stein manifold is a precise way of saying that the group of its 
automorphisms is big, the following problem of T\'oth and Varolin
\cite{TothVarolin2000} is very natural. 

%
%   THE TOTH-VAROLIN PROBLEM
%
\begin{problem}\label{prob:TothVarolin}
Is every Stein manifold with the density property which is diffeomorphic to $\R^{2n}$
also biholomorphic to $\C^n$? 
\end{problem}

This question remains unsolved. An affirmative answer would yield
(following Derksen and Kutzschebauch \cite{DerksenKutzschebauch1998}) 
a nonlinearizable holomorphic action of $\C^*$ on $\C^3$, and we would know that holomorphic 
$\C^*$-actions on $\C^n$ are linearizable if and only if $n\le 2$. 
As it stands, the case $n=3$ is unsolved, and it can only be solved 
with a good characterisation of $\C^3$ in hand. Incidentally, the linearization problem was 
the main motivation of Makar--Limanov for introducing his invariant. 
The most spectacular use of the Makar--Limanov invariant was the proof 
by Kaliman et al. \cite{Kalimanet.al.1997} that every algebraic 
$\C^*$-actions on $\C^3$ is linearizable. 
This invariant was the crucial tool in proving that all potential counterexamples to linearization 
(like the Koras--Russell 3-fold)  are non-isomorphic to $\C^3$.

It is still unknown whether an affine algebraic variety which is biholomorphic to $\C^n$ is also 
algebraically isomorphic to $\C^n$. This is known as  {\em Zaidenberg's problem} \cite{Zaidenberg}.
If Problem \ref{prob:TothVarolin} of T\'oth and Varolin has an affirmative answer, then 
the Koras--Russell cubic KR is a counterexample to Zaidenberg's problem. 
Indeed, Leuenberger \cite{Leuenberger2016} showed that KR has the density property
(see Example (5) in Subsec.\ \ref{ss:DP}), and it is known to be diffeomorphic to $\R^6$.

One can naturally generalise the cancellation problem as follows. 

\begin{problem}\label{prob:CP2}
Let $X$ and $Y$ be affine algebraic manifolds such that $X\times \C$ is algebraically isomorphic 
to $Y\times \C$. Does it follow that $X$ is isomorphic to $Y$? 
The analogous problem for Stein manifolds and biholomorphisms.
\end{problem}

Here the answer is negative in general, and additional conditions must be imposed.
For example, Danielewski found that for a polynomial $p$ with 
simple roots the affine surfaces $D_n := \{ (x,y,z)\in \C^3 : x^n y = p(x)\}$ satisfy 
$D_n \times \C \cong D_m \times \C$ (unpublished preprint, 1989). Later, Fieseler \cite{Fieseler1994} 
showed that $D_n$ and $D_m$ for $n\ne m$ are not even homeomorphic by examining their
fundamental group at infinity. 
This gives counterexamples to both the holomorphic and the algebraic cancellation. 
Cancellation also fails in the differentiable category. 
Counterexamples include some nice smooth complex algebraic varieties.
Take for example a surface of Ramanujam which is contractible but not simply connected 
at infinity. The surface 
\[
	 S:=\Big\{(x,y,z)\in \C^3 : \frac{ (xz+1)^3 - (yz+1)^2 -z }{z} = 0  \Big\} 
\]
is such an example. It is not homeomorphic to $\R^4$,
but $S\times \C$ is contractible and simply connected at infinity. 
This is easily seen from the definition of the fundamental group at infinity; 
however, we remark that in general for a smooth affine algebraic 
variety of dimension $n\ge 3$, contractibility implies simple connectedness at infinity \cite{Zaidenberg}.  
Thus, $S\times \R^2$ is diffeomorphic to $\R^6$. By Miyanishi--Sugie 
\cite{MiyanishiSugie1980} it follows that $S\times \C$ is not algebraically isomorphic to $\C^3$. 
It remains unknown whether $S\times \C$ is biholomorphic to $\C^3$.

The affine modification $M$ of $\C^4$ with divisor $D:=\C^3\times \{0\}\subset \C^4$ along the center 
$C:= S\subset \C^3\times \{0\}$ is another interesting example. It is given by one equation in $\C^5$:
\[
	M = \Big\{ (x,y,z,u,v) \in \C^5 : uv=\frac{(xz+1)^3 - (yz+1)^2 -z}{z}   \Big\}.
\]
The manifold $M$ is known to be diffeomorphic to $\R^8$ by works of Kaliman and Zaidenberg 
\cite{KalimanZaidenberg1999TG} and Kaliman and Kutzschebauch \cite{KalimanKutzschebauch2008IM}. 
If $M$ is isomorphic to $\C^4$, we have an algebraic $\C^4$ in $\C^5$
which is not algebraically straightenable since the defining polynomial 
for $M$ has singular fibres. No such examples are known. 
Except for the classical case of Abhyankar--Moh--Suzuki of lines in the plane 
(see the discussion after Corollary \ref{cor:nonstraightenable}), it is unknown whether a smooth 
algebraic hypersurface in $\C^{n+1}$ which is isomorphic to $\C^n$ is 
algebraically straightenable. If $M$ is not isomorphic to
$\C^4$ but is biholomorphic to $\C^4$, it provides a counterexample to Zaidenberg's problem. 
Finally, if $M$ is not biholomorphic to $\C^4$ then it is a counterexample 
to Problem  \ref{prob:TothVarolin} of T\'oth and Varolin. It is still unknown whether $M$ is isomorphic 
or biholomorphic to $\C^4$, which shows how difficult these questions are.  

Let us finally mention some characterisations of affine spaces. Unfortunately, they are not useful for 
solving any of the  problems we have formulated. 

It is classically known in both categories that the algebra of regular or holomorphic functions,
respectively, determines the underlying space. For the affine algebraic case this holds by 
definition (the space is the spectrum of the algebra of regular functions), 
and for reduced and irreducible Stein spaces it is a classical result attributed to 
Remmert in his habilitation thesis; see Grauert and Remmert 
\cite[Theorem 6, p.\ 184]{GrauertRemmert2004}. 
% In both cases, this is an application of Cartan--Serre Theorems A and B. 
If one considers only  the ring structure of the space of  holomorphic functions, 
one must also impose the condition that the ring homomorphism maps $\imath=\sqrt{-1}$ 
to itself (and not to $-\imath$); see the appendix of the paper \cite{Issa1966}
which was published by Heisuke Hironaka under the pseudonym Hei Iss'sa.  
This result is known as Bers's Theorem after 
Lipman Bers who first proved it for  open Riemann surfaces. 
The fact that the $\C$-algebra of meromorphic functions determines the
underlying Stein space is known as Iss'sa's Theorem \cite{Issa1966}. 
In the same paper it is shown that a field isomorphism of the spaces of meromorphic functions 
has to map $\imath$ to $\imath$ in order to induce a biholomorphism 
(and not an antibiholomorphism) between the Stein spaces. 

In the algebraic setting, there are no applicable characterisations other than those 
by Miyanishi--Sugie and Ramanujam mentioned above.
A recent result of Cantat, Regeta and Xie  \cite{Cantat2019families} 
states that if $X$ is a reduced connected affine algebraic variety 
whose automorphism group $\Autalg (X)$ is isomorphic to 
$\Autalg (\C^n)$ as an abstract group, then $X$ is isomorphic to $\C^n$. 
Under the stronger assumption that 
the automorphism groups are isomorphic as ind-groups (a notion introduced by Shafarevich in
\cite{Shafarevich1966,Shafarevich1981}), this conclusion was obtained 
earlier by Kraft \cite{Kraft2017}. Another result of Andrist and Kraft \cite{AndristKraft2014} 
says that if the semigroups of self-maps $\End (X)$ and $\End (\C^n$) are isomorphic 
then, up to an automorphism of the base field $\C$, $X$ is isomorphic to $\C^n$. 
(The result in \cite{AndristKraft2014} is proved over an algebraically closed field of arbitrary characteristic. 
The algebraic case over the complex numbers is already implicitly contained in \cite{Andrist2011}.)
In the holomorphic category, Andrist proved in \cite{Andrist2011} that if the semigroup of self-maps 
$\End(X)$ of a complex space is abstractly isomorphic to the semigroup $\End(\C^n)$ 
then $X$ is biholomorphic or antibiholomorphic to $\C^n$.  
The difference in the results over the complex numbers is that, in the holomorphic case, 
the field automorphisms are continuous.

%
% Comment by Rafael Andrist regarding the previous paragraph. It has been implemented.
%
\begin{comment}
The result in [19] \cite{AndristKraft2014} is proved over an algebraically closed field of arbitrary characteristic. 
The algebraic case over the complex numbers is already implicitly contained in [15] \cite{Andrist2011}. 
The difference in the results over the complex numbers is that in the holomorphic case, 
the field automorphisms are continuous.
\end{comment}

A result of Isaev and Kruzhilin  \cite{IsaevKruzhilin}, which compares to Kraft's result with the ind-group structure, 
says that if  the automorphism group of a complex manifold $X$ of dimension $n$ is isomorphic as a topological 
group to $\Aut(\C^n)$ then $X$ is biholomorphic to $\C^n$. The isomorphism of topological groups 
implies by classical results that there is a real analytic action of the unitary affine group $\R^{2n} \rtimes U_n$ 
on $X$ by holomorphic automorphisms, and classifying manifolds with such isometries is not too difficult.

Another not very applicable criterion comes from the work of Kutzschebauch, L\'arusson and Schwarz 
\cite{KutzschebauchLarussonSchwarz2017TG}. 
As in the previous result, it assumes some group action on $X$. 
If a Stein manifold $X$ carries the action of a reductive group $G$ so that the categorical quotient is 
stratified biholomorphic to the quotient of a large linear action of $G$ on $\C^n$, then $X$ is 
biholomorphic to $\C^n$, and there is a biholomorphism linearizing the action. 
The same is true without the assumption of largeness for finite groups, 
as well as groups $G$ having connected component $\C^*$ \cite{KutzschebauchSchwarz2021} 
or $G=\SL_2 (\C)$ \cite{KutzschebauchLarussonSchwarz2017TG}.

A characterisation  in the holomorphic case which is clearly worth to be explored further is due to 
Boc Thaler \cite[Theorem IV.15]{BocThaler2016}.

\begin{theorem} \label{th:BT}
Let $X$ be a Stein manifold with the density property. Then $X$ is
biholomorphic to $\C^n$ if and only if $X$ can be exhausted by Runge images of the ball.
\end{theorem}

\smallskip
\noindent {\bf Acknowledgement.} 
The authors wish to thank Rafael Andrist, Shulim Kaliman, Finnur L\'arusson, and Riccardo Ugolini 
for their helpful remarks and suggestions.

%%%%%%%%%%
%%%%%%%%%%
%%%%%%%%%%
%%%%%%%%%%   THE BIBLIOGRAPHY
%%%%%%%%%%
%%%%%%%%%%

%{\bibliographystyle{abbrv} \bibliography{references}}
%\begin{comment}

% \end{comment}

\end{document}